\documentclass[10pt]{article}

\usepackage{amsmath}
\usepackage{amsfonts}
\usepackage{amssymb}
\usepackage[english]{babel}
\usepackage{amsthm}

\def\pf{\par\noindent {\bf Proof}~\par\noindent}
\def\qed{~\hfill{$\square$}\pagebreak[1]\par\medskip\par}

\newcommand{\mR}{\mathbb{R}}
\newcommand{\mC}{\mathbb{C}}
\newcommand{\mN}{\mathbb{N}}

\newcommand{\p}{\partial}
\newcommand{\ux}{\underline{x}}
\newcommand{\uu}{\underline{u}}
\newcommand{\uom}{\underline{\omega}}
\newcommand{\pux}{\underline{\partial}}

\newcommand{\mcA}{\mathcal{A}}
\newcommand{\mcB}{\mathcal{B}}
\newcommand{\mcC}{\mathcal{C}}
\newcommand{\mcT}{\mathcal{T}}
\newcommand{\mcU}{\mathcal{U}}
\newcommand{\enb}{\overline{e_0}}
\newcommand{\Dbar}{\overline{D}}
\newcommand{\pxn}{\partial_{x_0}}
\newcommand{\onehalf}{\frac{1}{2}}
\newcommand{\la}{\langle}
\newcommand{\ra}{\rangle}

\newtheorem{lemma}{Lemma}
\newtheorem{proposition}{Proposition}
\newtheorem{remark}{Remark}

\numberwithin{equation}{section}
\numberwithin{theorem}{section}
\numberwithin{proposition}{section}
\numberwithin{lemma}{section}
\numberwithin{definition}{section}
\numberwithin{remark}{section}
\numberwithin{corollary}{section}

\begin{document}

\title{On a Chain of Harmonic and Monogenic Potentials in Euclidean Half--space}

\author{F.\ Brackx, H.\ De Bie, H.\ De Schepper}

\date{\small{Clifford Research Group, Department of Mathematical Analysis,\\ Faculty of Engineering and Architecture, Ghent University\\
Building S22, Galglaan 2, B-9000 Gent, Belgium\\}}

\maketitle

\begin{abstract}
\noindent In the framework of Clifford analysis, a chain of harmonic and monogenic potentials is constructed in the upper half of Euclidean space $\mR^{m+1}$, including a higher dimensional generalization of the complex logarithmic function. Their distributional limits at the boundary $\mR^{m}$ turn out to be well-known distributions such as the Dirac distribution, the Hilbert kernel, the fundamental solution of the Laplace and Dirac operators, the square root of the negative Laplace operator, and the like. It is shown how each of those potentials may be recovered from an adjacent kernel in the chain by an appropriate convolution with such a distributional limit. 
\end{abstract}

\maketitle


\section{Introduction}


Consider in the upper half $\mC_+ = \{ z = x+iy \in \mC : y > 0 \}$ of the complex plane, the logarithmic function 
$$
\ln z = \ln |z| + i \arg z, \qquad \mbox{Im} \, z > 0
$$
with
$$
\arg z = \frac{\pi}{2} - \arctan \frac{x}{y},  \quad y>0
$$
It is quite an interesting function; let us have a closer look at its properties.
\begin{itemize}
\item[(i)] The function $\ln z$ is holomorphic in $\mC_+$, i.e.\ it is a null solution of the Cauchy--Riemann operator 
$$
D = \frac{1}{2} \left ( \p_x + i \p_y \right )
$$
\item[(ii)] Its real and imaginary parts are conjugate harmonic functions in $\mC_+$ and the real part $\ln |z| = \ln r = \ln \sqrt{x^2+y^2}$ is, up to a constant, the fundamental solution of the two--dimensional Laplace operator
$$
\Delta_2 = \p_{xx}^2 + \p_{yy}^2 = 4 D \overline{D}
$$
where $\overline{D} = \frac{1}{2} \left ( \p_x - i \p_y \right )$ is the complex conjugate Cauchy--Riemann operator; in fact, we have, in distributional sense,
$$
\Delta_2 \left ( \frac{1}{2\pi} \ln |z| \right ) = \delta(z)
$$
with $\delta(z)$ the Dirac or delta distribution in $\mC$.
\item[(iii)] As a holomorphic function, $\ln z$ has a complex derivative in $\mC_+$, given by
$$
\frac{d}{dz} \ln z = \frac{1}{z}
$$
meaning that $\ln z$ is a holomorphic primitive (or potential) in $\mC_+$, with respect to the complex derivative $\frac{d}{dz}$, of the function $\frac{1}{z}$ which, in its turn, is, up to a constant, the fundamental solution of the Cauchy--Riemann operator; in fact we have, in distributional sense,
$$
D \left (\frac{1}{\pi} \frac{1}{z} \right ) = \delta(z)
$$
Note that the complex derivative operator $\frac{d}{dz}$ is nothing else but the conjugate Cauchy--Riemann operator $\overline{D}$, and there also holds in $\mC_+$
$$
\frac{d}{dz} \ln z = \overline{D} \ln z = \p_x \ln z = (-i \p_y) \ln z = \frac{1}{z}
$$
\item[(iv)] The conjugate harmonic real and imaginary parts $\ln r$ and $\arg z$ satisfy the Cauchy--Riemann system
$$
\left \{ \begin{array}{lllllll}
\p_x \ln r & = & \phantom{-} \p_y \arg z & = & \frac{x}{x^2+y^2} & = & \phantom{-} \mbox{Re} \left ( \frac{1}{z} \right ) \\[3mm]
\p_y \ln r & = & - \p_x \arg z & = & \frac{y}{x^2+y^2} & = & - \mbox{Im} \left ( \frac{1}{z} \right ) \\
\end{array} \right .
$$
where at the right hand sides one recognizes the Poisson kernel $P(x,y) = \frac{y}{x^2+y^2}$ and its harmonic conjugate 
$Q(x,y) = \frac{x}{x^2+y^2}$ in $\mC_+$; it follows that
$$
\overline{D} ( 2 \ln r) = \frac{1}{z}
$$
and 
$$
\overline{D} ( 2 i \arg z )  = \frac{1}{z}
$$
meaning that the functions $2 \ln r$ and $2i \arg z$ are conjugate harmonic potentials in $\mC_+$, with respect to the operator $\overline{D} = \frac{d}{dz}$, of the Cauchy kernel $\frac{1}{z}$.
\item[(v)] The distributional limits for $y \rightarrow 0+$ of the Cauchy kernel $\frac{1}{z}$ and its holomorphic potential $\ln z$, are given by
$$
\lim_{y \rightarrow 0+} \frac{1}{z} = \lim_{y \rightarrow 0+} \frac{x}{x^2+y^2} - i \lim_{y \rightarrow 0+} \frac{y}{x^2+y^2} = \mbox{Pv}\, \frac{1}{x} - i \pi \delta(x)
$$
with Pv$\frac{1}{x}$ the ''principal value'' distribution on the real axis, 
and 
$$
\lim_{y \rightarrow 0+} \ln z = \lim_{y \rightarrow 0+} \ln |z| + i \lim_{y \rightarrow 0+} \arg z = \ln |x| + i \pi Y(-x)
$$
with $Y(x)$ the Heaviside step function. These distributional boundary values fit into the following two commutative schemes
$$
\begin{array}{ccc}
\ln r & \xrightarrow{\hspace*{2mm} \p_x \hspace*{2mm}} & \frac{x}{x^2+y^2} \\[2mm]
\hspace*{-7mm} ^{\hspace*{0.1mm}_{y \rightarrow 0+}} \downarrow & & \downarrow \\
\ln |x| & \xrightarrow{\hspace*{2mm} \p_x \hspace*{2mm}} & \mbox{Pv} \, \frac{1}{x}
\end{array}
\qquad \mbox{and} \qquad  
\begin{array}{ccc}
- \arg z & \xrightarrow{\hspace*{2mm} \p_x \hspace*{2mm}} & \frac{y}{x^2+y^2} \\[2mm]
\hspace*{-5mm} ^{\hspace*{0.1mm}_{y \rightarrow 0+}} \downarrow & & \downarrow \\
- \pi Y(-x) & \xrightarrow{\hspace*{2mm} \p_x \hspace*{2mm}} & \pi \, \delta(x)
\end{array}
$$
and moreover they form Hilbert pairs, the Hilbert transform on the real axis being given by 
$$
\mathcal{H}[T] = H \ast T = \frac{1}{\pi} \mbox{Pv} \, \frac{1}{x} \ast T,
$$
since we have indeed 
$$
\mathcal{H} \left[ \pi \delta(x)\right] = \mbox{Pv} \, \frac{1}{x}, \qquad \mathcal{H} \left [ \mbox{Pv} \, \frac{1}{x} \right ] = \pi \delta(x) 
$$
and
$$
\mathcal{H} \left[ \ln |x| \right] = \pi \, Y(-x), \qquad \mathcal{H} \left[ \pi \, Y(-x) \right] = \ln |x|
$$
\end{itemize}

The aim of this paper is to construct a generalization of this logarithmic potential function in higher dimension, more specifically in the framework of Clifford analysis, where the functions under consideration take their values in the universal Clifford algebra $\mR_{0,m+1}$ constructed over Euclidean space $\mR^{m+1}$ equipped with a quadratic form of signature $(0,m+1)$. The concept of a higher dimensional holomorphic function, mostly called monogenic function, is expressed by means of a generalized Cauchy--Riemann operator, which is a combination of the derivative with respect to one of the variables, say $x_0$, and the so--called Dirac operator $\pux$ in the remaining variables $(x_1, x_2, \ldots, x_m)$. The generalized Cauchy--Riemann operator and its Clifford algebra conjugate linearize the Laplace operator, whence Clifford analysis is entitled to be qualified as a refinement of harmonic analysis.\\[-2mm]

It is a remarkable fact that the thus constructed monogenic logarithmic function in upper half--space $\mR^{m+1}_+$ shows the same, above mentioned, five properties as in the complex plane. Starting point of our construction is the fundamental solution of the generalized Cauchy--Riemann operator, also called Cauchy kernel, and its relation to the Poisson kernel and its harmonic conjugate in $\mR^{m+1}_+$. We then proceed by induction in two directions, {\em downstream} by differentiation and {\em upstream} by primitivation, yielding an doubly infinite chain of monogenic, and thus harmonic, potentials. This chain mimics the well--known sequence of holomorphic potentials in $\mC_+$ (see e.g. \cite{slang}):
$$
\frac{1}{k!} z^k \left[ \ln z - ( 1 + \frac{1}{2} + \ldots + \frac{1}{k}) \right] \rightarrow \ldots \rightarrow z ( \ln z - 1) \rightarrow \ln z 
\stackrel{\frac{d}{dz}}{\longrightarrow} \frac{1}{z} \rightarrow - \frac{1}{z^2} \rightarrow \ldots \rightarrow (-1)^{k-1} \frac{(k-1)!}{z^k}
$$
Identifying the boundary of upper half--space with $\mR^m \cong \{(x_0,\ux) \in \mR^{m+1} : x_0 = 0\}$, the distributional limits for $x_0 \rightarrow 0+$ of those potentials are computed; they divide into two classes which are linked by the Hilbert transform and encompass well--known distributions in $\mR^m$ such as the Dirac or delta distribution, the Hilbert kernel, the fundamental solutions of the Dirac and the Laplace operators, the square root of the negative Laplacian, and the like. It is also shown how each of the monogenic potentials may be recovered from an adjacent kernel in the chain by an appropriate convolution with such a boundary distribution.\\[-2mm]

The organization of the paper is as follows. To make the paper self--contained we recall in Section 2 the basics of Clifford algebra and Clifford analysis. In Section 3 we construct a conjugate harmonic in upper half--space $\mR^{m+1}_+$ to the fundamental solution of the $(m+1)$--dimensional Laplace operator, which is essential to obtaining the desired monogenic logarithmic function in $\mR^{m+1}_+$. In Section 4 we study the so--called downstream potentials obtained under the action of the Clifford algebra conjugate of the generalized Cauchy--Riemann operator. Finally, in Section 5, we study the monogenic logarithmic function in $\mR^{m+1}_+$ and we construct, by an appropriate form of primitivation, the sequence of upstream potentials. Section 6 is concluding.


\section{Basics of Clifford analysis}


Clifford analysis (see e.g. \cite{red}) is a function theory which offers a natural and elegant generalization to higher dimension of holomorphic functions in the complex plane and refines harmonic analysis. Let  $(e_0, e_1,\ldots,e_m)$ be the canonical orthonormal basis of Euclidean space $\mR^{m+1}$ equipped with a quadratic form of signature $(0,m+1)$. Then the non--commutative multiplication in the universal real Clifford algebra $\mR_{0,m+1}$ is governed by the rule 
$$
e_{\alpha} e_{\beta} + e_{\beta} e_{\alpha} = -2 \delta_{\alpha \beta}, \qquad \alpha,\beta = 0, 1,\ldots,m
$$
whence $\mR_{0,m+1}$ is generated additively by the elements $e_A = e_{j_1} \ldots e_{j_h}$, where $A=\lbrace j_1,\ldots,j_h \rbrace \subset \lbrace 0,\ldots,m \rbrace$, with $0\leq j_1<j_2<\cdots < j_h \leq m$, and $e_{\emptyset}=1$. 
For an account on Clifford algebra we refer to e.g. \cite{porteous}.\\[-2mm]

We identify the point $(x_0, x_1, \ldots, x_m) \in \mR^{m+1}$ with the Clifford--vector variable 
$$
x = x_0 e_0  + x_1 e_1  + \cdots x_m e_m = x_0 e_0  + \ux
$$ 
and the point $(x_1, \ldots, x_m) \in \mR^{m}$ with the Clifford--vector variable $\ux$. 
Introducing spherical co--ordinates $\ux = r \uom$, $r = |\ux|$, $\uom \in S^{m-1}$, gives rise to the Clifford--vector valued locally integrable function $\uom$, which is to be seen as the higher dimensional analogue of the {\em signum}--distribution on the real line; we will encounter $\uom$ as one of the distributions discussed below.\\[-2mm]

At the heart of Clifford analysis lies the so--called Dirac operator 
$$
\p = \pxn e_0 + \p_{x_1} e_1 + \cdots \p_{x_m} e_m =  \pxn e_0 + \pux
$$
which squares to the negative Laplace operator: $\p^2 = - \Delta_{m+1}$, while also $\pux^2 = - \Delta_{m}$. Due to the non--commutative character of the multiplication in the Clifford algebra, the Dirac operator may act from the left or from the right on a Clifford algebra valued function with, in general, different results. The (left and right) fundamental solution of the Dirac operator $\p$ is given by
$$
E_{m+1} (x) = - \frac{1}{\sigma_{m+1}} \ \frac{x}{|x|^{m+1}}
$$
where $\sigma_{m+1} = \frac{2\pi^{\frac{m+1}{2}}}{\Gamma(\frac{m+1}{2})}$ stands for the area of the unit sphere $S^{m}$ in $\mR^{m+1}$.
We also introduce the generalized Cauchy--Riemann operator 
$$
D = \onehalf \enb \p = \onehalf (\pxn + \enb \pux)
$$ 
and its Clifford algebra conjugate $\Dbar = \onehalf(\pxn - \enb \pux)$. As is the case in the complex plane, both operators decompose the Laplace operator in $\mR^{m+1}$: $D \Dbar = \Dbar D = \frac{1}{4} \Delta_{m+1}$. \\[-2mm]

A continuously differentiable function $F(x)$, defined in an open region $\Omega \subset \mR^{m+1}$ and taking its values in the Clifford algebra  $\mR_{0,m+1}$, is called (left--)monogenic if it satisfies the equation $D F = 0$ in $\Omega$, which is equivalent with $\p F = 0$. \\

Singling out the basis vector $e_0$, we can decompose the real Clifford algebra $\mR_{0,m+1}$ in terms of the Clifford algebra 
$\mR_{0,m}$ as $\mR_{0,m+1} = \mR_{0,m} \oplus \enb \mR_{0,m}$. Similarly we decompose the functions considered as 
$$
F(x_0,\ux) = F_1(x_0,\ux) + \enb F_2(x_0,\ux)
$$ 
where $F_1$ and $F_2$ take their values in the Clifford algebra $\mR_{0,m}$; mimicking functions of a complex variable, we will call $F_1$ the {\em real} part and $F_2$ the {\em imaginary} part of the function $F$.\\[-2mm]

We will extensively use two families of distributions in $\mR^m$, which have been thoroughly studied in \cite{fb1,fb2,distrib}. The first family $\mathcal{T} = \{ T_\lambda : \lambda \in \mC\}$ is very classical. It consists of the radial distributions 
$$
T_\lambda = {\rm Fp} \ r^{\lambda} = {\rm Fp} \ (x_1^2 + \ldots + x_m^2)^{\frac{\lambda}{2}}
$$
their action on a test function $\phi \in \mathcal{S}(\mR^m)$ being given by
$$
\la T_\lambda, \phi \ra = \sigma_m \la {\rm Fp} \; r^\mu_+, \Sigma^{(0)}[\phi] \ra 
$$
with $\mu = \lambda +m-1$. In the above expressions ${\rm Fp}\; r^\mu_+$ is the classical {\em finite part} distribution on the real $r$-axis  and $\Sigma^{(0)}$ is the scalar valued generalized spherical mean, defined on scalar valued test functions $\phi(\ux)$ by
$$
\Sigma^{(0)}[\phi] = \frac{1}{\sigma_m} \int_{S^{m-1}} \phi(\ux) \, dS(\uom)
$$
This family $\mathcal{T}$ contains a.o. the fundamental solution of the Laplace operator. As convolution operators they give rise to the traditional Riesz potentials (see e.g. \cite{helgason}). The second family $\mathcal{U} = \{ U_\lambda : \lambda \in \mC\}$ of distributions arises in a natural way by the action of the  Dirac operator $\pux$ on $\mathcal{T}$. The $U_{\lambda}$--distributions thus are typical Clifford analysis objects: they are Clifford--vector valued, and they also arise as products of $T_{\lambda}$--distributions with the distribution $\uom = \frac{\ux}{|\ux|}$, mentioned above. The action of $U_\lambda$ on a test function $\phi \in \mathcal{S}(\mR^m)$ is given by
$$
\la U_\lambda, \phi \ra  = \sigma_m \la {\rm Fp} \; r^\mu_+, \Sigma^{(1)}[\phi] \ra 
$$
with $\mu = \lambda +m-1$, and where the Clifford--vector valued generalized spherical mean $\Sigma^{(1)}$ is defined on scalar valued test functions $\phi(\ux)$ by
$$
\Sigma^{(1)}[\phi] = \frac{1}{\sigma_m} \int_{S^{m-1}} \uom \ \phi(\ux) \, dS(\uom) 
$$
Typical example in the $\mathcal{U}$--family is the fundamental solution of the Dirac operator.\\[-2mm]

The normalized distributions $T^{*}_\lambda$ and $U^{*}_\lambda$ arise when the singularities of $T_\lambda$ and $U_\lambda$ are removed by dividing them by an appropriate Gamma-function, showing the same simple poles. The $T^{*}_\lambda$--distributions are defined by
\begin{eqnarray*}
\left \{
\begin{array}{ll}
\displaystyle{T_\lambda^* = \pi^{\frac{\lambda+m}{2}} \frac{T_\lambda}{\Gamma \left ( \frac{\lambda+m}{2} \right )}}, & \lambda \ne -m-2l\\[5mm]
\displaystyle{T_{-m-2l}^* = \frac{\pi^{\frac{m}{2}-l}}{2^{2l} \Gamma \left ( \frac{m}{2} + l \right )} (-\Delta)^l \delta (\ux)}, & l \in \mN_0
\end{array}
\right . 
\end{eqnarray*}
while the Clifford--vector valued distributions $U^{*}_\lambda$ are defined by
\begin{eqnarray*}
\left \{
\begin{array}{ll}
\displaystyle{U_\lambda^* = \pi^{\frac{\lambda+m+1}{2}} \, \frac{U_\lambda}{\Gamma \left ( \frac{\lambda + m + 1}{2} \right )}}, & \lambda \ne -m-2l-1\\[5mm]
\displaystyle{U_{-m-2l-1}^* = - \frac{\pi^{\frac{m}{2}-l}}{2^{2l+1} \, \Gamma \left ( \frac{m}{2} + l + 1 \right )} \; \pux^{2l+1} \delta(\ux)}, & l \in \mN_0
\end{array}
\right . 
\end{eqnarray*}

The normalized distributions $T_\lambda^*$ and $U_\lambda^*$ are holomorphic mappings from $\lambda \in \mC$ to the space $\mathcal{S}'(\mR^m)$ of tempered distributions. As already mentioned they are intertwined by the action of the Dirac operator. They enjoy the following properties: for all $\lambda \in \mC$ one has
\begin{itemize}
\item[(i)] $\ux \; T_\lambda^* = \frac{\lambda+m}{2\pi} \; U_{\lambda+1}^*$; \quad
$\ux \; U_\lambda^* = U_\lambda^* \; \ux = - T_{\lambda+1}^*$
\item[(ii)] $\pux \; T_\lambda^* = \lambda \; U_{\lambda-1}^*$; \quad 
$\pux \; U_\lambda^* = U_\lambda^* \; \pux = - 2\pi \; T_{\lambda-1}^*$
\item [(iii)] $\Delta_m T_\lambda^* = 2 \pi \lambda T_{\lambda-2}^*$ ; \quad
$\Delta_m U_\lambda^* = 2 \pi (\lambda-1) U_{\lambda-2}^*$
\item[(iv)] $r^2 T_\lambda^* = \frac{\lambda + m}{2\pi} \ T_{\lambda+2}^*$; \quad
$r^2 U_\lambda^* = \frac{\lambda + m + 1}{2\pi} \ U_{\lambda+2}^*$
\end{itemize}

Of particular importance for the sequel are the convolution formulae for the $T_\lambda^*$-- and $U_\lambda^*$--distributions; we list them in the following proposition and refer the reader to \cite{distrib} for more details.

\begin{proposition}
\label{prop1}
\rule{0mm}{0mm}
\begin{itemize}
\item[(i)] For all $(\alpha,\beta) \in \mC \times \mC$ such that $ \alpha \neq 2j, j \in \mN_0$, $ \beta \neq 2k, k \in \mN_0$ and $\alpha+\beta+m \neq 2l, l \in \mN_0$ the convolution $T_{\alpha}^* \ast T_{\beta}^*$ is the tempered distribution given by
$$
T_{\alpha}^* \ast T_{\beta}^* = c_m(\alpha,\beta)\; T_{\alpha+\beta+m}^*
$$
with
\begin{eqnarray*}
c_m(\alpha,\beta) & = & \pi^{\frac{m}{2}}\; \frac{\Gamma \left( -\frac{\alpha + \beta + m}{2} \right)}{\Gamma \left( -\frac{\alpha}{2} \right) \Gamma \left( -\frac{\beta}{2} \right)}
\end{eqnarray*}

\item[(ii)] For $(\alpha,\beta) \in \mC \times \mC$ such that $ \alpha \neq 2j+1, \  \beta \neq 2k, \ \alpha+\beta \neq -m+2l+1, \ j, k, l \in \mN_0$ one has
$$
U_{\alpha}^*  \ast   T_{\beta}^*  =   T_{\beta}^* \ast  U_{\alpha}^* =  c_m(\alpha - 1,\beta)\; U_{\alpha+\beta+m}^*
$$

\item[(iii)] For $(\alpha,\beta) \in \mC \times \mC$ such that $ \alpha \neq 2j+1, \  \beta \neq 2k+1, \ \alpha+\beta \neq -m+2l, \ j, k, l \in \mN_0$ one has
$$
U_{\alpha}^*  \ast   U_{\beta}^*	  = U_{\beta}^*  \ast   U_{\alpha}^*	  =  \pi^{\frac{m}{2}+1} \displaystyle{ \frac{\Gamma(- \frac{\alpha+\beta+m}{2})}{\Gamma(\frac{-\alpha+1}{2}) \Gamma(\frac{-\beta+1}{2})}}  \;  T_{\alpha+\beta+m}^*
$$
\end{itemize}
\end{proposition}

\begin{remark}
The action of a Clifford algebra valued distribution on a ditto test function is assumed to be carried out componentwise, the respective basis vectors being multiplied in the Clifford algebra.
\end{remark}

\begin{remark}
In general the convolution of Clifford algebra valued distributions is not commutative. However, as is seen from formula (iii) in Proposition \ref{prop1}, the convolution of two distributions from the $\mcU$--family is indeed commutative. We will frequently use this property in the sequel. Convolution by distributions from the $\mcT$--family is intrinsically commutative since they are scalar valued.
\end{remark}

\begin{remark}
In general the convolution of distributions is not associative. However, as is seen from the formulae in Proposition \ref{prop1}, the convolution of distributions from the  $\mcT$-- and $\mcU$--families is associative. Also this property will be frequently used this in the sequel. 
\end{remark}


\section{A conjugate harmonic to Green's function}
\label{conjharm}


The fundamental solution of the Laplace operator $\Delta_{m+1}$ in $\mR^{m+1}$, sometimes called Green's function, and here denoted, for reasons which  will become clear afterwards, by $\frac{1}{2}A_0(x_0,\ux)$, is given by
\begin{equation}
\frac{1}{2}A_0(x_0,\ux) = - \frac{1}{m-1} \frac{1}{\sigma_{m+1}} \frac{1}{|x|^{m-1}}
\label{A0}
\end{equation}
Considering the function $A_0(x_0,\ux)$ as a harmonic function in the upper half--space $\mR^{m+1}_+$, our aim now is to construct its conjugate harmonic in $\mR^{m+1}_+$ in the sense of \cite{red}, in this way elaborating further on an earlier result of \cite{Xu}. This means that we have to look for a harmonic function $B_0(x_0,\ux)$ in $\mR^{m+1}_+$ such that
$$
C_0(x_0,\ux) = \frac{1}{2} A_0(x_0,\ux) + \frac{1}{2} \overline{e_0} \, B_0(x_0,\ux)
$$
is monogenic in $\mR^{m+1}_+$ w.r.t.\ the generalized Cauchy--Riemann operator $D$. Expressing the monogenicity of $C_0$ in $\mR^{m+1}_+$ leads to the system
\begin{equation}
\left \{ \begin{array}{rcl}
\p_{x_0} A_0 + \pux B_0 & = & 0 \\[2mm]
\p_{x_0} B_0 + \pux A_0 & = & 0 
\end{array} \right .
\label{CR12}
\end{equation}
which clearly mimics the Cauchy--Riemann system in the complex plane. Taking into account the explicit expression (\ref{A0}) of $A_0(x_0,\ux)$, the system (\ref{CR12}) reduces to
\begin{equation}
\left \{ \begin{array}{rcl}
\pux B_0(x_0,\ux) &=& - \displaystyle\frac{2}{\sigma_{m+1}} \, \displaystyle\frac{x_0}{|x|^{m+1}} \ = \ - P(x_0,\ux)  \\[3mm] 
\p_{x_0} B_0(x_0,\ux) &=& - \displaystyle\frac{2}{\sigma_{m+1}} \, \displaystyle\frac{\ux}{|x|^{m+1}} \ = \ \phantom{-} Q(x_0,\ux) 
\end{array} \right .
\label{CR34}
\end{equation}
where $P$ and $Q$ stand for the Poisson kernel and its conjugate in $\mR^{m+1}_+$ (see also Section \ref{downstream}). From the second condition in (\ref{CR34}) its follows that, for an arbitrary, but fixed, $x_0^\ast$,
$$
B(x_0,\ux) = \frac{2}{\sigma_{m+1}} \, \frac{\ux}{|\ux|^m} \, F_m \left ( \frac{|\ux|}{x_0} \right ) - \frac{2}{\sigma_{m+1}} \, \frac{\ux}{|\ux|^m} \, F_m \left ( \frac{|\ux|}{x_0^\ast} \right ) + W(\ux)
$$
where we have put
$$
F_m(v) = \int_0^v \frac{\eta^{m-1}}{(1+\eta^2)^\frac{m+1}{2}} \, d\eta = \frac{v^m}{m} \,  _2F_1 \left ( \frac{m}{2},\frac{m+1}{2};\frac{m}{2}+1;-v^2 \right )
$$
with $_2F_1$ a standard hypergeometric function (see e.g. \cite{grad}).\\[-2mm]

From the first condition in (\ref{CR34}) it then follows that the function $W(\ux)$ should satisfy the equation
$$
\pux W(\ux) = - \left ( \p_{x_0} A_0 \right )_{x_0^\ast} = - \frac{2}{\sigma_{m+1}} \, \frac{x_0^\ast}{|x_0^\ast e_0 + \ux|^{m+1}}
$$
and a straightforward calculation shows that the function
$$
W(\ux) = \frac{2}{\sigma_{m+1}} \, \frac{\ux}{|\ux|^m} \, F_m \left ( \frac{|\ux|}{x_0^\ast} \right )
$$
does the job. A conjugate harmonic to $A_0$ in $\mR_+^{m+1}$ is thus given by
\begin{equation}
B_0(x_0,\ux) = \frac{2}{\sigma_{m+1}} \, \frac{\ux}{|\ux|^m} \, F_m \left ( \frac{|\ux|}{x_0} \right )
\label{B0}
\end{equation}
or
\begin{equation}
B_0(x_0,\ux) = \frac{2}{m} \, \frac{1}{\sigma_{m+1}} \, \frac{\ux}{x_0^m} \, _2F_1 \left ( \frac{m}{2},\frac{m+1}{2};\frac{m}{2}+1;-\frac{|\ux|^2}{x_0^2} \right )
\label{B02}
\end{equation}
Expression (\ref{B02}) clearly shows that $B_0(x_0,\ux)$ is well--defined for $\ux=0$, with
$$
\lim_{\ux \rightarrow 0} B_0(x_0,\ux) = 0, \qquad x_0 > 0
$$
Taking into account that
$$
F_m(+\infty) = \int_0^{+\infty} \frac{\eta^{m-1}}{(1+\eta^2)^\frac{m+1}{2}} \, d\eta = \frac{\sqrt{\pi}}{2} \frac{\Gamma \left ( \frac{m}{2} \right )}{\Gamma \left ( \frac{m+1}{2} \right )}
$$
expression (\ref{B0}) leads to the following distributional limit
\begin{equation}
b_0(\ux) = \lim_{x_0 \rightarrow 0+} B_0(x_0,\ux) = \frac{1}{\sigma_m} \frac{\ux}{|\ux|^m} = \frac{1}{\pi} \, \frac{1}{\sigma_m} \, U^\ast_{-m+1}
\label{b0}
\end{equation}
in which one recognizes, up to a minus sign, the fundamental solution $E_m(\ux)$ of the Dirac operator $\pux$ in $\mR^m$:
$$
- b_0(\ux) = \frac{1}{\sigma_m} \, \frac{\overline{\ux}}{|\ux|^m} = - \frac{1}{\pi} \, \frac{1}{\sigma_m} \, U^\ast_{-m+1} = E_m(\ux)
$$
This distribution $E_m(\ux)$ may act as a convolution kernel for the so--called $T$--operator, which is a convolution operator acting on Clifford algebra valued Schwartz--functions $f \in \mathcal{S}(\mR^m)$ or on ditto tempered distributions as
$$
T[f] = E_m \ast f = - b_0 \ast f
$$
Seen the fact that $E_m(\ux)$ is the fundamental solution of the Dirac operator $\pux$, this $T$--operator is an inverse to this Dirac operator:
$$
\pux \, T[f] = f, \qquad f \in \mathcal{S}(\mR^m)
$$
The Green function $A_0(x_0,\ux)$ itself shows the following distributional limit:
\begin{equation}
a_0(\ux) = \lim_{x_0 \rightarrow 0+} A_0(x_0,\ux) = - \frac{2}{m-1} \, \frac{1}{\sigma_{m+1}} \, {\rm Fp} \frac{1}{|\ux|^{m-1}} = 
- \frac{2}{m-1} \, \frac{1}{\sigma_{m+1}} \, T_{-m+1}^\ast
\label{a0}
\end{equation}
Using this distribution, up to a minus sign, as a convolution kernel, gives rise to the convolution operator $(-\Delta)^{-\frac{1}{2}}$, acting on Schwartz--functions or tempered distributions by, see e.g.\ \cite{helgason},
$$
\left ( - \Delta \right )^{- \frac{1}{2}} [f] = \frac{2}{m-1} \, \frac{1}{\sigma_{m+1}} \, T^\ast_{-m+1} \ast f = - a_0 \ast f
$$
The functions $A_0(x_0,\ux)$ and $B(x_0,\ux)$ being conjugate harmonic in $\mR^{m+1}_{+}$, we expect their distributional boundary values $a_0(\ux)$ and $b_0(\ux)$ to be intimately related. This is indeed the case, as will be shown in Section \ref{upstream}.


\section{Downstream potentials}
\label{downstream}

\subsection{The Cauchy kernel as a potential}

As is well--known, the Cauchy kernel of Clifford analysis, i.e. the fundamental solution of the generalized Cauchy--Riemann operator $D$,
$$
C_{-1}(x_0,\ux) = \frac{1}{\sigma_{m+1}} \, \frac{x \overline{e_0}}{|x|^{m+1}} =  \frac{1}{\sigma_{m+1}} \, \frac{x_0 - 
\overline{e_0} \ux}{|x|^{m+1}}
$$
may be decomposed in terms of the Poisson kernels in $\mR^{m+1}_+$:
$$
C_{-1}(x_0,\ux) = \frac{1}{2} A_{-1}(x_0,\ux) + \frac{1}{2} \overline{e_0} \, B_{-1}(x_0,\ux)
$$
where, also mentioning the traditional notations, for $x_0 >0$,

\begin{equation}
\left \{ \begin{array}{rcl}
A_{-1}(x_0,\ux) & = & P(x_0,\ux) \ = \ \phantom{-} \frac{2}{\sigma_{m+1}} \, \frac{x_0}{|x|^{m+1}}\\[4mm]
B_{-1}(x_0,\ux) & = & Q(x_0,\ux) \ = \ - \frac{2}{\sigma_{m+1}} \, \frac{\ux}{|x|^{m+1}}
\end{array} \right .
\label{A-1B-1}
\end{equation}

Note that the Poisson kernel $A_{-1}$ is real--valued, while its conjugate harmonic kernel $B_{-1}$ is Clifford vector--valued. Their distributional limits for $x_0 \rightarrow 0+$ are given by
\begin{eqnarray*}
a_{-1}(\ux) & = & \lim_{x_0 \rightarrow 0+} A_{-1}(x_0,\ux) \ = \ \delta(\ux) \ \ = \ \phantom{-} \frac{2}{\sigma_m} \, T^\ast_{-m} \\
b_{-1}(\ux) & = & \lim_{x_0 \rightarrow 0+} B_{-1}(x_0,\ux) \ = \ H(\ux) \ = \ - \frac{2}{\sigma_{m+1}} \, U^\ast_{-m} 
\end{eqnarray*}
and also
$$
c_{-1}(\ux) = \lim_{x_0 \rightarrow 0+} C_{-1}(x_0,\ux) \ = \ \frac{1}{2} \delta(\ux) + \frac{1}{2} \overline{e_0} \, H(\ux)
$$
Note that the distribution
$$
H(\ux) = - \frac{2}{\sigma_{m+1}} \, U^\ast_{-m} = - \frac{2}{\sigma_{m+1}} \, \mbox{Pv} \frac{\ux}{|\ux|^{m+1}}
$$
where $\mbox{Pv}$ stands for the {\em principal value} distribution in $\mR^m$, is the convolution kernel of the Hilbert transform $\mathcal{H}$ in $\mR^m$ (see e.g.\ \cite{gilmur}). Note also that both distributional boundary values are linked by this Hilbert transform:
\begin{eqnarray*}
\mathcal{H} \left [ a_{-1} \right ] & = & \mathcal{H} \left [ \delta \right ] \ = \ H \ast \delta \ = \ H \ = \ b_{-1} \\
\mathcal{H} \left [ b_{-1} \right ] & = & \mathcal{H} \left [ H \right ] \ = \ H \ast H \ = \ \delta \ = \ a_{-1}
\end{eqnarray*}
since $\mathcal{H}^2 = \mathbf{1}$, while$$\overline{e_0} \, \mathcal{H} \left [ c_{-1} \right ] = c_{-1}$$
Conversely, the Poisson kernels are the Poisson transforms of these distributional limits:
\begin{eqnarray*}
\mathcal{P} \left [ a_{-1} \right ] & = & P(x_0,\cdot) \ast a_{-1}(\cdot)(\ux) \ = \ P(x_0,\cdot) \ast \delta(\cdot)(\ux) \ = \ P(x_0,\ux)\\
\mathcal{P} \left [ b_{-1} \right ] & = & P(x_0,\cdot) \ast b_{-1}(\cdot)(\ux) \ = \ P(x_0,\cdot) \ast H(\cdot)(\ux) \ = \ Q(x_0,\ux)
\end{eqnarray*}
It follows that also the Poisson kernels themselves are linked by the Hilbert transform in the variable $\ux \in \mR^m$:
\begin{eqnarray*}
\mathcal{H} \left [ A_{-1} \right ] & = & H(\cdot) \ast A_{-1}(x_0,\cdot)(\ux) \ = \ H(\cdot) \ast P(x_0,\cdot)(\ux) \\
& = & P(x_0,\cdot) \ast H(\cdot)(\ux) \ = \ Q(x_0,\ux) \ = \ B_{-1}(x_0,\ux) \\[2mm]
\mathcal{H} \left [ B_{-1} \right ] & = & \mathcal{H}^2 \left [ A_{-1} \right ] \ = \ A_{-1}
\end{eqnarray*}
For a function $f \in L_2(\mR^m)$, its Poisson transforms
\begin{eqnarray*}
\mathcal{P}[f] & = & P(x_0,\cdot) \ast f(\cdot)(\ux) \ =\ A_{-1}(x_0,\cdot) \ast f(\cdot)(\ux) \\
\mathcal{Q}[f] & = & Q(x_0,\cdot) \ast f(\cdot)(\ux) \ = \ B_{-1}(x_0,\cdot) \ast f(\cdot)(\ux)
\end{eqnarray*}
belong to the Clifford--Hardy space $\mbox{Harm}^2 \left ( \mR^{m+1}_+ \right )$ of Clifford algebra valued harmonic functions in $\mR_+^{m+1}$:
$$
\mbox{Harm}^2 \left( \mR^{m+1}_+ \right) = \left \{ F(x_0,\ux) \, : \, F \mbox{\ is harmonic in $\mR_+^{m+1}$ and\ } \sup_{x_0 > 0} \int_{\mR^m} |F(x_0,\ux)|^2 \, d\ux \; < \; + \infty \right \}
$$
and show the non--tangential $L_2$--boundary values $$\lim_{x_0 \rightarrow 0+} \mathcal{P}[f] = f \quad {\rm and} \quad \lim_{x_0 \rightarrow 0+} \mathcal{Q}[f] = \mathcal{H}[f]$$ with
$$
\mathcal{H}[f] = H \ast f = \frac{2}{\sigma_{m+1}} \, \mbox{Pv} \int_{\mR^m} \frac{\overline{\uu}}{|\uu|^{m+1}} \, f(\ux - \uu) \, d\uu
$$
the explicit expression for the Hilbert transform of $f$. In its turn the Cauchy transform of $f \in L_2(\mR^m)$, given by
$$
\mathcal{C}[f] = C_{-1}(x_0,\cdot) \ast f(\cdot)(\ux) = \frac{1}{2} \mathcal{P}[f] + \frac{1}{2} \overline{e_0} \mathcal{Q}[f]
$$
belongs to the Clifford--Hardy space $H^2(\mR^{m+1}_+)$ of monogenic functions in $\mR^{m+1}_+$:
$$
H^2 \left ( \mR^{m+1}_+ \right ) = \left \{ F(x_0,\ux) \, : \, F \mbox{\ is monogenic in $\mR_+^{m+1}$ and\ } \sup_{x_0 > 0} \int_{\mR^m} |F(x_0,\ux)|^2 \, d\ux \; < \; + \infty \right \}
$$
and shows the following non--tangential $L_2$--boundary value:
$$
\lim_{x_0 \rightarrow 0+} \mathcal{C}[f] = \frac{1}{2} f + \frac{1}{2} \overline{e_0} \mathcal{H}[f] = \left ( \frac{1}{2} \delta + \frac{1}{2} \overline{e_0} H \right ) \ast f = \mathcal{A}\mathcal{S}[f]
$$
which belongs to the Clifford--Hardy space $H^2(\mR^m)$, see \cite{gilmur}. In signal analysis the functions in $H^2(\mR^m)$ are called {\em analytic signals} ; they show no negative-frequency components (see e.g. \cite{hahn}). Whence the notation $\mathcal{A}\mathcal{S}$ for the boundary value of the Cauchy transform. Note that for this Cauchy transform we have several equivalent expressions:
$$
\mathcal{C} [f] = \mathcal{C} \left  [\overline{e_0} \, \mathcal{H}[f] \right ] = \mathcal{C} \left [ \mathcal{A}\mathcal{S}[f] \right ] = \mathcal{P} \left [ \mathcal{A}\mathcal{S}[f] \right ]
$$
From the monogenicity of the Cauchy kernel $C_{-1}(x_0,\ux)$ in $\mR^{m+1}_+$, i.e. 
$$
DC_{-1} = \frac{1}{2} \left ( \p_{x_0} + \overline{e_0} \pux \right ) C_{-1}=0
$$ 
it follows that the Poisson kernels $A_{-1}(x_0,\ux)$ and $B_{-1}(x_0,\ux)$ satisfy the generalized Cauchy--Riemann system
\begin{equation}
\left \{ \begin{array}{rcl}
\p_{x_0} A_{-1} + \pux B_{-1} & = & 0 \\[2mm]
\p_{x_0} B_{-1} + \pux A_{-1} & = & 0 \\
\end{array} \right .
\label{GCR}
\end{equation}
and that
\begin{equation}
\overline{D} C_{-1} = \frac{1}{2} \p_{x_0} C_{-1} - \frac{1}{2} \overline{e_0} \pux C_{-1} = \p_{x_0} C_{-1} = - \overline{e_0} \pux C_{-1}
\label{DC-1}
\end{equation}
and also that
\begin{equation}
\left \{ \begin{array}{l}
\overline{D} A_{-1} = \frac{1}{2} \p_{x_0} A_{-1} - \frac{1}{2} \overline{e_0} \pux A_{-1} = \frac{1}{2} \p_{x_0} A_{-1} + \frac{1}{2} \overline{e_0} \p_{x_0} B_{-1} = \p_{x_0} C_{-1} = \overline{D} C_{-1} \\[2mm]
\overline{D}(\overline{e_0} B_{-1}) = \frac{1}{2} \overline{e_0} \p_{x_0} B_{-1} - \frac{1}{2} \overline{e_0} \pux \overline{e_0} B_{-1} = \frac{1}{2} \overline{e_0} \p_{x_0} B_{-1} + \frac{1}{2} \p_{x_0} A_{-1} = \p_{x_0} C_{-1} = \overline{D} C_{-1}
\end{array} \right .
\label{DAB-1}
\end{equation}
Now we put
$$
\overline{D} C_{-1} = C_{-2} = \frac{1}{2} A_{-2} + \frac{1}{2} \overline{e_0} B_{-2}
$$
clearly a monogenic function in $\mR^{m+1}_+$, since $DC_{-2} = D \overline{D} C_{-1} = \frac{1}{4} \Delta_{m+1} C_{-1} = 0$. From this definition it follows that
$$
\left \{ \begin{array}{l}
A_{-2} = \p_{x_0} A_{-1} = - \pux B_{-1} \\[2mm]
B_{-2} = \p_{x_0} B_{-1} = - \pux A_{-1}
\end{array} \right .
$$
leading to the explicit expressions for the conjugate harmonic components of $C_{-2}$:
\begin{equation}
\left \{ \begin{array}{rcl}
A_{-2} & = & \displaystyle\frac{2}{\sigma_{m+1}} \, \displaystyle\frac{1}{|x|^{m+3}} \left ( |x|^2 - (m+1) x_0^2 \right ) = \displaystyle\frac{2}{\sigma_{m+1}} \, \displaystyle\frac{- m x_0^2 + |\ux|^2}{|x|^{m+3}} \\[4mm]
B_{-2} & = & (m+1) \, \displaystyle\frac{2}{\sigma_{m+1}} \, \displaystyle\frac{x_0 \ux}{|x|^{m+3}}
\end{array} \right .
\label{A-2B-2}
\end{equation}
Note that $A_{-2}(x_0, \ux)$ is real--valued, while $B_{-2}(x_0, \ux)$ is Clifford vector--valued. Moreover it is readily confirmed that they satisfy the generalized CR--system 
$$
\left \{ \begin{array}{l}
\p_{x_0} A_{-2} + \pux B_{-2} = 0 \\[2mm]
\p_{x_0} B_{-2} + \pux A_{-2} = 0 
\end{array} \right .
$$
The above relations (\ref{DC-1})--(\ref{DAB-1}) imply that the monogenic function $C_{-2}(x_0,\ux)$ in $\mR^{m+1}_+$ shows the monogenic potential (or primitive) $C_{-1}(x_0,\ux)$ and the conjugate harmonic potentials $A_{-2}(x_0,\ux)$ and $\overline{e_0} B_{-2}(x_0,\ux)$. The distributional limits for $x_0 \rightarrow 0+$ of these harmonic potentials are given by
$$
\left \{ \begin{array}{rcl}
a_{-2}(\ux) = \lim_{x_0 \rightarrow 0+} A_{-2}(x_0,\ux) & = & \displaystyle\frac{2}{\sigma_{m+1}} \, {\rm Fp} \displaystyle\frac{1}{|\ux|^{m+1}} \ = \ - \displaystyle\frac{4 \pi}{\sigma_{m+1}} T^\ast_{-m-1}\\[4mm]
b_{-2}(\ux) = \lim_{x_0 \rightarrow 0+} B_{-2}(x_0,\ux) & = & - \pux \delta \ = \ \displaystyle\frac{2m}{\sigma_m} \, U^\ast_{-m-1}
\end{array} \right .
$$
Conversely, the harmonic potentials $A_{-2}(x_0,\ux)$ and $B_{-2}(x_0,\ux)$ are recovered from these distributional boundary values by the Poisson transform 
\begin{eqnarray*}
A_{-2}(x_0,\ux) & = & \mathcal{P} \left [ a_{-2}(\ux) \right ] \ = \ P(x_0,\cdot) \ast a_{-2}(\cdot)(\ux) \ = \ A_{-1}(x_0,\cdot) \ast a_{-2}(\cdot)(\ux) \ = \ a_{-2}\ast A_{-1} \\
& = & \mathcal{Q} \left [ b_{-2}(\ux) \right ] \ = \ Q(x_0,\cdot) \ast b_{-2}(\cdot)(\ux) \ = \ B_{-1}(x_0,\cdot) \ast b_{-2}(\cdot)(\ux) \ = \ b_{-2} \ast B_{-1}
\end{eqnarray*}
and
\begin{eqnarray*}
B_{-2}(x_0,\ux) & = & \mathcal{P} \left [ b_{-2}(\ux) \right ] \ = \ P(x_0,\cdot) \ast b_{-2}(\cdot)(\ux) \ = \ A_{-1}(x_0,\cdot) \ast b_{-2}(\cdot)(\ux) \ = \ b_{-2} \ast A_{-1}\\
& = & \mathcal{Q} \left [ a_{-2}(\ux) \right ] \ = \ Q(x_0,\cdot) \ast a_{-2}(\cdot)(\ux) \ = \ B_{-1}(x_0,\cdot) \ast a_{-2}(\cdot)(\ux) \ = \ a_{-2} \ast B_{-1} 
\end{eqnarray*}
In the distribution $a_{-2}$ one recognizes the convolution kernel $-\pux H = - H \pux$, known as the Hilbert--Dirac kernel, see \cite{fbhds}, or perhaps better known as the convolution kernel for the pseudodifferential operator $(-\Delta)^\frac{1}{2}$ (see \cite{helgason}). The distribution $b_{-2}$ is, up to a minus sign, the Dirac derivative of the delta-distribution. Both distributional boundary values are linked by the Hilbert transform, as shown a.o.\ in the following lemma.

\begin{lemma}
One has
\begin{itemize}
\item[(i)] $- \pux a_{-1} = b_{-2}$, $-\pux b_{-1} = a_{-2}$, $- \overline{e_0} \pux c_{-1} = c_{-2}$
\item[(ii)] $\mathcal{H} \left [ a_{-2} \right ] = b_{-2}$, $\mathcal{H} \left [ b_{-2} \right ] = a_{-2}$, $\overline{e_0} \mathcal{H} \left [ c_{-2} \right ] = c_{-2}$
\item[(iii)] $c_{-1} \ast a_{-2} = c_{-2}$, $c_{-1} \ast \overline{e_0} b_{-2} = c_{-2}$, $c_{-1} \ast c_{-2} = c_{-2}$
\end{itemize}
\end{lemma}

\pf

\noindent (i) Follows by direct calculation.\\[-2mm]

\noindent (ii) Making use of the convolution calculation rules, recalled in Proposition \ref{prop1}, we have
\begin{eqnarray*}
\mathcal{H} \left [ a_{-2} \right ] & = & - \frac{4 \pi}{\sigma_{m+1}} \, H \ast T^\ast_{-m-1} \ = \ \frac{8 \pi}{\sigma_{m+1}^2} U^\ast_{-m} \ast T^\ast_{-m-1} \\
& = & \frac{8 \pi}{\sigma_{m+1}^2} \, \pi^\frac{m}{2} \, \frac{\Gamma \left ( \frac{m+2}{2} \right )}{\left ( \Gamma \left ( \frac{m+1}{2} \right ) \right )^2} \, U^\ast_{-m-1} \ = \ \frac{2m}{\sigma_m} \, U^\ast_{-m-1} \ = \ - \pux \delta \ = \ b_{-2}
\end{eqnarray*}
and
\begin{eqnarray*}
\mathcal{H} \left [b_{-2} \right ] & = & \frac{2m}{\sigma_m} \, H \ast U^\ast_{-m-1} \ = \ - \frac{4m}{\sigma_m \sigma_{m+1}} U^\ast_{-m} \ast U^\ast_{-m-1} \\
& = & - \frac{4m}{\sigma_m \sigma_{m+1}} \, \pi^\frac{m}{2} \, \frac{\Gamma \left ( \frac{m+1}{2} \right )}{\Gamma \left ( \frac{m+1}{2} \right ) \Gamma \left ( \frac{m+2}{2} \right )} \, T^\ast_{-m-1} \ = \ - \frac{4}{\sigma_{m+1}} \, T^\ast_{-m-1} \ = \ a_{-2}
\end{eqnarray*}

\noindent (iii) We subsequently find
$$
c_{-1} \ast a_{-2} = \left ( \frac{1}{2} \delta + \frac{1}{2} \overline{e_0} H \right ) \ast a_{-2} = \frac{1}{2} a_{-2} + \frac{1}{2} \overline{e_0} \mathcal{H} \left [ a_{-2} \right ] = \frac{1}{2} a_{-2} + \frac{1}{2} \overline{e_0} b_{-2} = c_{-2}
$$
$$
c_{-1} \ast \overline{e_0} b_{-2} = \left ( \frac{1}{2} \delta + \frac{1}{2} \overline{e_0} H \right ) \ast \overline{e_0} b_{-2} = \frac{1}{2} \overline{e_0} b_{-2} + \frac{1}{2} \mathcal{H} \left [ b_{-2} \right ] = \frac{1}{2} \overline{e_0} b_{-2} + \frac{1}{2} a_{-2} = c_{-2}
$$
and
$$
c_{-1} \ast c_{-2} = c_{-1} \ast \left ( \frac{1}{2} a_{-2} + \frac{1}{2} \overline{e_0} b_{-2} \right ) = \frac{1}{2} c_{-2} + \frac{1}{2} c_{-2} = c_{-2}
$$
\qed

Through the Poisson transform, the Hilbert--link between the distributional boundary values $a_{-2}$ and $b_{-2}$ is reflected in a similar relationship between the harmonic potentials $A_{-2}$ and $B_{-2}$, as it was also the case for $A_{-1}$ and $B_{-1}$. Indeed, one has
\begin{eqnarray*}
\mathcal{H} \left [ A_{-2} \right ] & = & H(\cdot) \ast A_{-2}(x_0,\cdot)(\ux) \ = \ H(\cdot) \ast P(x_0,\cdot) \ast a_{-2}(\cdot)(\ux) \\
& = & P(x_0,\cdot) \ast H(\cdot) \ast a_{-2}(\cdot)(\ux) \ = \ P(x_0,\cdot) \ast b_{-2}(\cdot)(\ux) \ = \ B_{-2}(x_0,\ux)\\
\mathcal{H} \left [ B_{-2} \right ] & = & \mathcal{H}^2 \left [ A_{-2} \right ] = A_{-2}
\end{eqnarray*}
These relations may also be rewritten as
\begin{eqnarray*}
b_{-1}(\cdot) \ast A_{-2}(x_0,\cdot)(\ux) = B_{-2}(x_0,\ux) \\
b_{-1}(\cdot) \ast B_{-2}(x_0,\cdot)(\ux) = A_{-2}(x_0,\ux)
\end{eqnarray*}
while, quite trivially,
\begin{eqnarray*}
a_{-1}(\cdot) \ast A_{-2}(x_0,\cdot)(\ux) = A_{-2}(x_0,\ux) \\
a_{-1}(\cdot) \ast B_{-2}(x_0,\cdot)(\ux) = B_{-2}(x_0,\ux)
\end{eqnarray*}
Note the following two commutative schemes, which are each others Hilbert image:
\begin{equation}
\begin{array}{ccc}
A_{-1} & \xrightarrow{\hspace*{2mm} -\pux \hspace*{2mm}} & B_{-2} \\[2mm]
\hspace*{-8mm} ^{\hspace*{0.1mm}_{x_0 \rightarrow 0+}} \downarrow & & \downarrow \\
\delta = a_{-1} & \xrightarrow{\hspace*{2mm} -\pux \hspace*{2mm}} & b_{-2} = - \pux \delta
\end{array}
\qquad \qquad \mbox{and} \qquad  \qquad 
\begin{array}{ccc}
B_{-1} & \xrightarrow{\hspace*{2mm} -\pux \hspace*{2mm}} & A_{-2} \\[2mm]
\hspace*{-7mm} ^{\hspace*{0.1mm}_{x_0 \rightarrow 0+}} \downarrow & & \downarrow \\
H = b_{-1} & \xrightarrow{\hspace*{2mm} -\pux \hspace*{2mm}} & a_{-2} = - \pux H
\end{array}
\label{cscheme1}
\end{equation}
By means of the distributional limits $a_{-2}$ and $b_{-2}$, we are now able to prove some remarkable relations between the conjugate harmonic components of $C_{-1}$ and $C_{-2}$; in fact they are shown to be linked by the distributional limits $a_0$ and $b_0$ of the Green function and its conjugate (see Section \ref{conjharm}).
\begin{proposition}
\label{propdownstream}
One has, convolutions being taken in the variable $\ux \in \mR^m$:
\begin{itemize}
\item[(i)] $a_0(\cdot) \ast A_{-2}(x_0,\cdot)(\ux) = A_{-1}(x_0,\ux) = b_0(\cdot) \ast B_{-2}(x_0,\ux)$
\item[(ii)]  $a_0(\cdot) \ast B_{-2}(x_0,\cdot)(\ux) = B_{-1}(x_0,\ux) = b_0(\cdot) \ast A_{-2}(x_0,\ux)$
\item[(iii)]  $a_0(\cdot) \ast C_{-2}(x_0,\cdot)(\ux) = C_{-1}(x_0,\ux) = b_0(\cdot) \ast C_{-2}(x_0,\ux)$
\item[(iv)]  $c_0(\cdot) \ast A_{-2}(x_0,\cdot)(\ux) = c_0(\cdot) \ast B_{-2}(x_0,\ux) = c_0(\cdot) \ast C_{-2}(x_0,\ux) = C_{-1}(x_0,\ux)$
\end{itemize}
\end{proposition}

\pf
\noindent (i)(ii) Put $b_0 \ast B_{-2} = A'_{-1}$. Then
\begin{eqnarray*}
\p_{x_0} A'_{-1} & = & b_0 \ast \p_{x_0} B_{-2} = b_0 \ast \left ( - \pux A_{-2} \right ) = - b_0 \pux \ast A_{-2} = \delta \ast A_{-2} = A_{-2} \\
\pux A'_{-1} & = & \pux b_0 \ast B_{-2} = - \delta \ast B_{-2} = B_{-2}
\end{eqnarray*}
while moreover
$$
\lim_{x_0 \rightarrow 0+} A'_{-1} = b_0 \ast b_{-2} = b_0 \ast (-\pux \delta) = - b_0 \pux \ast \delta = \delta \ast \delta = \delta = a_{-1}
$$
Similarly, by putting $b_0 \ast A_{-2} = B'_{-1}$, we have
\begin{eqnarray*}
\p_{x_0} B'_{-1} & = & b_0 \ast \p_{x_0} A_{-2} = b_0 \ast \left ( - \pux B_{-2} \right ) = - b_0 \pux \ast B_{-2} = \delta \ast B_{-2} = B_{-2} \\
\pux B'_{-1} & = & \pux b_0 \ast A_{-2} = - \delta \ast A_{-2} = A_{-2}
\end{eqnarray*}
while moreover
$$
\lim_{x_0 \rightarrow 0+} B'_{-1} = b_0 \ast a_{-2} = b_0 \ast (-\pux H) = - b_0 \pux \ast H = \delta \ast H = H = b_{-1}
$$
So $A'_{-1}$ and $B'_{-1}$ satisfy the CR--system (\ref{GCR}) and show the same distributional limits for $x_0 \rightarrow 0+$ as $A_{-1}$ and $B_{-1}$, respectively, from which it follows that they have to coincide: $A'_{-1} = A_{-1}$ and $B'_{-1} = B_{-1}$. Now note that
\begin{eqnarray*}
a_0 \ast A_{-2} & = & a_0 \ast \mathcal{H} \left [ B_{-2} \right ] = a_0 \ast H \ast B_{-2} = H \ast a_0 \ast B_{-2} = b_0 \ast B_{-2} = A_{-1} \\
a_0 \ast B_{-2} & = & a_0 \ast \mathcal{H} \left [ A_{-2} \right ] = a_0 \ast H \ast A_{-2} = H \ast a_0 \ast A_{-2} = b_0 \ast A_{-2} = B_{-1}
\end{eqnarray*}
to complete the proof of (i) and (ii).\\[-2mm]

\noindent (iii)(iv) It suffices to make the appropriate combinations of the results in (i) and (ii).
\qed

For a function $f \in L_2(\mR^m)$ we can define in $\mR^{m+1}_+$ the conjugate harmonic functions 
$$
\mathcal{A}_{-2} [f] = A_{-2}(x_0,\cdot) \ast f(\cdot) (\ux) \qquad \mbox{and} \qquad \mathcal{B}_{-2}[f] = B_{-2}(x_0,\cdot) \ast f(\cdot) (\ux)
$$
and the monogenic function
$$
\mathcal{C}_{-2} [f]  = C_{-2}(x_0,\cdot) \ast f(\cdot) (\ux)
$$
They show non--tangential $L_2$--boundary values on condition that $f$ belongs to the Clifford--Sobolev space 
$$
W_2^1(\mR^m) = \left \{ f \in L_2(\mR^m) \, : \, \pux f \in L_2(\mR^m) \right \}
$$ 
Under these assumptions there holds
\begin{eqnarray*}
\mcA_{-2}^+[f] = \lim_{x_0 \rightarrow 0+} \mathcal{A}_{-2}[f] &=& a_{-2} \ast f \ = \ - H \pux \ast f \ = \ - \mathcal{H} \left [ \pux f \right ] \ = \ - \pux \mathcal{H}[f] \\
\mcB_{-2}^+[f]  = \lim_{x_0 \rightarrow 0+} \mathcal{B}_{-2}[f] &=& b_{-2} \ast f \ = \ - \pux \delta \ast f \ = \ - \pux f 
\end{eqnarray*}
and also
$$
\mcC_{-2}^+[f]  = \lim_{x_0 \rightarrow 0+} \mathcal{C}_{-2}[f] = - \frac{1}{2} \pux \mathcal{H} [f] - \frac{1}{2} \overline{e_0} \pux f = (- \overline{e_0} \pux) \left ( \frac{1}{2} f + \frac{1}{2} \overline{e_0} \mathcal{H}[f] \right ) = (- \overline{e_0} \pux) \left ( \mathcal{A} \mathcal{S} [f] \right )
$$
Note that the convolution operators $\mathcal{A}_{-2}^+$,  $\mathcal{B}_{-2}^+$ and $\mathcal{C}_{-2}^+$ are bounded operators from $W_2^1(\mR^m)$ into $L_2(\mR^m)$, and that, for $f \in W_2^1(\mR^m)$, the $L_2$--boundary value 
$\mathcal{C}_{-2}^+[f]$ belongs to the Clifford--Hardy space $H^2(\mR^m)$. Also note the following commutative scheme for a function $f \in W_2^1(\mR^m)$:
$$
\begin{array}{rcl}
\mathcal{C}_{-1}[f] & \xrightarrow{\hspace*{2mm} -\overline{e_0} \pux \hspace*{2mm}} & \mathcal{C}_{-2}[f] \\[2mm]
^{\hspace*{0.1mm}_{x_0 \rightarrow 0+}} \downarrow \hspace*{4mm} & & \hspace*{4mm} \downarrow ^{\hspace*{0.1mm}_{x_0 \rightarrow 0+}}\\
\frac{1}{2} f + \frac{1}{2} \overline{e_0} \mathcal{H}[f] = \mathcal{A}\mathcal{S}[f] & \xrightarrow{\hspace*{2mm} -\overline{e_0} \pux \hspace*{2mm}} & - \overline{e_0} \pux \mathcal{A} \mathcal{S}[f] = - \frac{1}{2} \pux \mathcal{H}[f] - \frac{1}{2} \overline{e_0} \pux f
\end{array}
$$
which reflects at the level of the operators, the commutative schemes (\ref{cscheme1}) at the level of the convolution kernels.

\subsection{Further derived potentials}

Proceeding in the same manner as in Subsection 4.1, we can define a sequence of monogenic potentials in $\mR_+^{m+1}$:
$$
C_{-k-1} = \overline{D} C_{-k} = \overline{D}^2 C_{-k+1} = \ldots = \overline{D}^k C_{-1}, \qquad k=1,2,\ldots
$$
where each monogenic potential decomposes into two conjugate harmonic potentials:
$$
C_{-k-1} = \frac{1}{2} A_{-k-1} + \frac{1}{2} \overline{e_0} B_{-k-1}, \qquad k=1,2,\ldots
$$
with, for $k$ odd, say $k=2\ell-1$,
$$
\left \{ \begin{array}{rcl}
A_{-2\ell} & = & \p_{x_0}^{2\ell-1} A_{-1} \ = \ - \p_{x_0}^{2\ell-2} \pux B_{-1} \ = \ \ldots \ = \ - \pux^{2\ell-1} B_{-1} \\[2mm]
B_{-2\ell} & = & \p_{x_0}^{2\ell-1} B_{-1} \ = \ - \p_{x_0}^{2\ell-2} \pux A_{-1} \ = \ \ldots \ = \ - \pux^{2\ell-1} A_{-1} 
\end{array} \right .
$$
while for $k$ even, say $k=2\ell$,
$$
\left \{ \begin{array}{rcl}
A_{-2\ell-1} & = & \p_{x_0}^{2\ell} A_{-1} \ = \ - \p_{x_0}^{2\ell-1} \pux B_{-1} \ = \ \ldots \ = \ \pux^{2\ell} A_{-1} \\[2mm]
B_{-2\ell-1} & = & \p_{x_0}^{2\ell} B_{-1} \ = \ - \p_{x_0}^{2\ell-1} \pux A_{-1} \ = \ \ldots \ = \ \pux^{2\ell} B_{-1} 
\end{array} \right .
$$
Note that also holds
$$
C_{-k-1} = \p_{x_0} C_{-k} = (-\overline{e_0} \pux) C_{-k} =  \p_{x_0}^2 C_{-k+1} = (-\overline{e_0} \pux)^2 C_{-k+1} = \ldots =  \p_{x_0}^k C_{-1} = (-\overline{e_0} \pux)^k C_{-1}
$$
while the conjugate harmonic components satisfy the recurrence relations
\begin{equation}
\left \{ \begin{array}{rcl}
A_{-k-1} & = & \p_{x_0} A_{-k} = - \pux B_{-k} \\[2mm]
B_{-k-1} & = & \p_{x_0} B_{-k} = - \pux A_{-k} 
\end{array} \right .
\label{AB-k-1}
\end{equation}
whence
$$
\left \{ \begin{array}{rcl}
\overline{D} A_{-k} & = & \displaystyle\frac{1}{2} \left ( \p_{x_0} - \overline{e_0} \pux \right ) A_{-k} = \displaystyle\frac{1}{2} A_{-k-1} + \displaystyle\frac{1}{2} \overline{e_0} B_{-k-1} = C_{-k-1} \\[4mm]
\overline{D} ( \overline{e_0} B_{-k}) & = & \displaystyle\frac{1}{2} \left ( \p_{x_0} \overline{e_0} - \pux \right ) B_{-k} = \displaystyle\frac{1}{2} A_{-k-1} + \displaystyle\frac{1}{2} \overline{e_0} B_{-k-1} = C_{-k-1}
\end{array} \right .
$$
which expresses the fact that $A_{-k}$ and $\overline{e_0} B_{-k}$ are indeed potentials (or primitives) of $C_{-k-1}$. Their distributional limits for $x_0 \rightarrow 0+$ are given by
$$
\left \{ \begin{array}{rcl}
a_{-2\ell} & = & (- \pux)^{2\ell-1} H  = - 2^{2\ell-1} \displaystyle\frac{\Gamma \left ( \frac{m+2\ell-1}{2} \right )} {\pi^{\frac{m-2\ell+1}{2}}} \, T^{\ast}_{-m-2\ell+1}  \\[5mm]
& = & (-1)^{\ell-1}  2^{\ell-1} (2\ell-1)!!  \displaystyle\frac{\Gamma \left ( \frac{m+2\ell-1}{2} \right )} {\pi^{\frac{m+1}{2}}} \, {\rm Fp} \displaystyle\frac{1}{r^{m+2\ell-1}}                                \\[7mm]
b_{-2\ell} & = & (- \pux)^{2\ell-1} \delta = 2^{2\ell-1} \displaystyle\frac{\Gamma \left (\frac{m+2\ell}{2} \right )} {\pi^{\frac{m-2\ell+2}{2}}} \, U^{\ast}_{-m-2\ell+1}       \end{array} \right .
$$
and
$$
\left \{ \begin{array}{rcl}
a_{-2\ell-1} & = &  \pux^{2\ell} \delta = 2^{2\ell} \displaystyle\frac{\Gamma \left (\frac{m+2\ell}{2} \right )} {\pi^\frac{m-2\ell}{2}}  \, T^{\ast}_{-m-2\ell}                    
\\[7mm]
b_{-2\ell-1} & = & \pux^{2\ell} H  \ = \ - 2^{2\ell} \displaystyle\frac{\Gamma \left ( \frac{m+2\ell+1}{2} \right )} {\pi^\frac{m-2\ell+1}{2}} \, U^{\ast}_{-m-2\ell}  \\[5mm]
& = &(-1)^{\ell-1}  2^{\ell} (2\ell-1)!!  \displaystyle\frac{\Gamma \left ( \frac{m+2\ell+1}{2} \right )}{\pi^\frac{m+1}{2}} \,  {\rm Fp} \displaystyle\frac{1}{r^{m+2\ell}} \, \omega  \end{array} \right .
$$
They show the following properties, which can be verified by direct calculation.
\begin{lemma}
\label{lem2}
One has for $j,k=1,2,\ldots$
\begin{itemize}
\item[(i)] $a_{-k} \xrightarrow{\hspace*{1mm} -\pux \hspace*{1mm}} b_{-k-1} \xrightarrow{\hspace*{1mm} -\pux \hspace*{1mm}} a_{-k-2}$
\item[(ii)] $\mathcal{H} \left [ a_{-k} \right ] = b_{-k}$, $\mathcal{H} \left [ b_{-k} \right ] = a_{-k}$
\item[(iii)] $a_{-j} \ast a_{-k} = a_{-j-k+1}$ \\
$a_{-j} \ast b_{-k} = b_{-j} \ast a_{-k} = b_{-j-k+1}$ \\
$b_{-j} \ast b_{-k} = a_{-j-k+1}$.
\end{itemize}
\end{lemma}
Through the Poisson transform, the above Hilbert--link (Lemma \ref{lem2}(ii)) between the distributional boundary values $a_{-k}$ and $b_{-k}$ is reflected into a similar relationship between the harmonic potentials $A_{-k}$ and $B_{-k}$, as was already shown for $k=1$ and $k=2$. Indeed, we have, the Hilbert transform being taken in the variable $\ux \in \mR^m$:
$\mathcal{H} \left [ A_{-k} \right ] = B_{-k}$ and $\mathcal{H} \left [ B_{-k} \right ] = A_{-k}$, which may also be written as
\begin{equation}
\label{ter}
\left \{ \begin{array}{rcl}
b_{-1} (\cdot) \ast A_{-k}(x_0,\cdot)(\ux) & = & B_{-k}(x_0,\ux) \\[2mm]
b_{-1} (\cdot) \ast B_{-k}(x_0,\cdot)(\ux) & = & A_{-k}(x_0,\ux)
\end{array} \right .
\end{equation}
while, trivially,
\begin{equation}
\label{quater}
\left \{ \begin{array}{rcl}
a_{-1} (\cdot) \ast A_{-k}(x_0,\cdot)(\ux) & = & A_{-k}(x_0,\ux) \\[2mm]
a_{-1} (\cdot) \ast B_{-k}(x_0,\cdot)(\ux) & = & B_{-k}(x_0,\ux)
\end{array} \right .
\end{equation}
Note also the following commutative scheme:
\begin{equation}
\begin{array}{ccccc}
A_{-k} & \xrightarrow{\hspace*{2mm} -\pux \hspace*{2mm}} & B_{-k-1} & \xrightarrow{\hspace*{2mm} -\pux \hspace*{2mm}}  & A_{-k-2} \\[2mm]
\hspace*{-8mm} ^{\hspace*{0.1mm}_{x_0 \rightarrow 0+}} \downarrow & & \downarrow & & \downarrow \\
a_{-k} & \xrightarrow{\hspace*{2mm} -\pux \hspace*{2mm}} & b_{-k-1} & \xrightarrow{\hspace*{2mm} -\pux \hspace*{2mm}}  & a_{-k-2} 
\end{array}
\label{cscheme2}
\end{equation}
The formulae (\ref{ter}) and (\ref{quater}) are special cases of the more general, and remarkable, result that the distributional boundary values $a_{-k}$ and $b_{-k}$ may act as convolution operators to convert the harmonic potentials into harmonic potentials of a lower order.
\begin{proposition}
One has for $k=1,2,\ldots$ and $j=0,1,\ldots$
\begin{itemize}
\item[(i)] $b_{-j-1} \ast A_{-k} = B_{-k-j}$
\item[(ii)] $b_{-j-1} \ast B_{-k} = A_{-k-j}$
\item[(iii)] $a_{-j-1} \ast A_{-k} = A_{-k-j}$
\item[(iv)] $a_{-j-1} \ast B_{-k} = B_{-k-j}$
\end{itemize}
\end{proposition}

\pf

\noindent (i) First assume that $k$ is even, say $k=2\ell$, and put $b_{-j-1} \ast A_{-2\ell} = B'_{-2\ell-j}$. Then we have
$$
\p_{x_0} B'_{-2\ell-j} = b_{-j-1} \ast \p_{x_0} A_{-2\ell} = b_{-j-1} \ast \left ( -\pux B_{-2\ell} \right ) = - \left ( b_{-j-1} \pux \right ) \ast B_{-2 \ell} = a_{-j-2} \ast  B_{-2\ell} 
$$
Assuming $j$ to be even, say $j=2i$, there holds
$$
a_{-2i-2} \ast B_{-2\ell} = (-\pux)^{2i+1} H \ast B_{-2\ell} = (-\pux)^{2i+1} A_{-2\ell} = B_{-2\ell-2i-1} = B_{-2\ell-j-1}
$$
while for $j$ odd, say $j=2i-1$, we have
$$
a_{-2i-2} \ast B_{-2\ell} = \pux^{2i} \delta \ast B_{-2\ell} = (-\pux)^{2i} B_{-2\ell} = B_{-2\ell-2i} = B_{-2\ell-j-1}
$$
and so $\p_{x_0} B'_{-2\ell-j} = B_{-2\ell - j-1}$. For the action of the Dirac operator $\pux$ on $B'_{-2\ell-j}$ we obtain
$$
\pux B'_{-2\ell-j} = \pux b_{-j-1} \ast A_{-2\ell} = - a_{-j-2} \ast A_{-2\ell}
$$
where now for $j$ even, say $j=2i$,
$$
-a_{-2i-2} \ast A_{-2\ell} = - (-\pux)^{2i+1} H \ast A_{-2 \ell} = - (-\pux)^{2i+1} B_{-2 \ell} = - A_{-2 \ell -2i -1} = - A_{-2\ell -j -1}
$$
while for $j$ odd, say $j=2i-1$,
$$
-a_{-2i-1} \ast A_{-2\ell} = - \pux^{2i} \delta \ast A_{-2\ell} = - (-\pux)^{2i} A_{-2\ell} = A_{-2\ell-2i} = - A_{-2\ell-j-1}
$$
and hence $\pux B'_{-2\ell-j} = -A_{-2\ell-j-1}$. Moreover
$$
\lim_{x_0 \rightarrow 0+} B'_{-2\ell-j} = b_{-j-1} \ast a_{-2\ell} = b_{-j-1} \ast (- \pux)^{2\ell-1} H = b_{-j-1} (- \pux)^{2\ell-1} \ast H = a_{-2\ell - j} \ast H = b_{-2\ell-j}
$$
This means that $B'_{-2\ell-j}$ and $B_{-2\ell-j}$ satisfy the same system (\ref{AB-k-1}) and show the same distributional boundary value for $x_0 \rightarrow 0+$, so they have to coincide, whence $b_{-j-1} \ast A_{-2\ell} = B_{-2\ell-j}$.\\
Next assume that $k$ is odd, say $k=2\ell-1$, and put $b_{-j-1} \ast A_{-2\ell+1} = B''_{-2\ell-j+1}$. Then we have, with a similar calculation as above: $\p_{x_0} B''_{-2\ell-j+1} = B_{-2\ell -j}$ and $\pux B''_{-2\ell-j+1} = - A_{-2\ell-j}$, while $\lim_{x_0 \rightarrow 0+} B''_{-2\ell-j+1} = b_{-2\ell-j+1}$, from which it follows that indeed $b_{-j-1} \ast A_{-2\ell+1} = B_{-2\ell - j+1}$, which completes the proof of (i). \\[-2mm]

\noindent (ii) The proof of (ii) is similar to that of (i).\\[-2mm]

\noindent (iii) Using the Hilbert--link between the distributional boundary values $a_{-j-1}$ and $b_{-j-1}$, we have indeed
$$
a_{-j-1} \ast A_{-k} = \mathcal{H} \left [ b_{-j-1} \right ] \ast A_{-k} = H \ast b_{-j-1} \ast A_{-k} = b_{-j-1} \ast H \ast A_{-k} = b_{-j-1} \ast B_{-k} = A_{-k-j}
$$
Note however that an alternative proof of (iii) is already contained in the calculations in the proofs of (i) and (ii).\\[-2mm]

\noindent (iv) The proof of (iv) is similar to that of (iii), now making use of (ii).
\qed

For a function $f \in L_2(\mR^m)$ we can define in $\mR^{m+1}_+$ the conjugate harmonic functions
$$
\mathcal{A}_{-k}[f] = A_{-k}(x_0,\cdot) \ast f(\cdot) (\ux) \qquad \mbox{and} \qquad \mathcal{B}_{-k}[f] = B_{-k}(x_0,\ux) \ast f(\cdot) (\ux), \quad k=1,2,\ldots
$$
and the monogenic function
$$
\mathcal{C}_{-k}[f] = C_{-k}(x_0,\cdot) \ast f(\cdot) (\ux), \qquad k=1,2,\ldots
$$
By definition of the potential kernels $A_{-k}$, $B_{-k}$ and $C_{-k}$, it is readily obtained that $\mathcal{A}_{-k}[f]$ and $\mathcal{B}_{-k}[f]$ are conjugate harmonic potentials of $\mathcal{C}_{-k-1}[f]$, while $\mathcal{C}_{-k}[f]$ is a monogenic potential (or primitive) of $\mathcal{C}_{-k-1}[f]$ in $\mR^{m+1}_+$. These potentials will show non--tangential $L_2$--boundary values for $x_0 \rightarrow 0+$ on condition that $f$ belongs to the Clifford--Sobolev space
$$
W_2^{k-1}(\mR^m) = \left \{ f \in L_2(\mR^m) \, : \, \pux f, \; \pux^2 f, \; \ldots, \pux^{k-1} f \in L_2(\mR^m) \right \}
$$
Under this assumption we have
\begin{eqnarray*}
\mathcal{A}_{-2\ell}^+ [f] = \lim_{x_0 \rightarrow 0+} \mathcal{A}_{-2\ell} [f] & = & a_{-2\ell} \ast f = - \pux^{2\ell-1} \mathcal{H}[f] = - \mathcal{H} \left [ \pux^{2\ell-1} f  \right ] \\[2mm]
\mathcal{B}_{-2\ell} ^+ [f]  = \lim_{x_0 \rightarrow 0+} \mathcal{B}_{-2\ell} [f] & = & b_{-2\ell} \ast f = - \pux^{2\ell-1} f \\[2mm]
\mathcal{A}_{-2\ell-1}^+ [f] = \lim_{x_0 \rightarrow 0+} \mathcal{A}_{-2\ell-1} [f] & = & a_{-2\ell-1} \ast f = \pux^{2\ell} f \\[2mm]
\mathcal{B}_{-2\ell-1}^+ [f]  = \lim_{x_0 \rightarrow 0+} \mathcal{B}_{-2\ell-1} [f] & = & b_{-2\ell-1} \ast f = \pux^{2\ell} \mathcal{H}[f] = \mathcal{H} \left [ \pux^{2\ell} f  \right ]
\end{eqnarray*}
and the convolution operators $\mathcal{A}_{-k}^+$ and $\mathcal{B}_{-k}^+$ are bounded operators from $W_2^{k-1}(\mR^m)$ into $L_2(\mR^m)$. For a function $f \in W_2^{k-1}(\mR^m)$ we also obtain the following expressions of the non--tangential $L_2$--boundary values of the monogenic potentials:
\begin{eqnarray*}
\mathcal{C}_{-2\ell}^+ [f]  =  \lim_{x_0 \rightarrow 0+} \mathcal{C}_{-2\ell} [f] & = & - \frac{1}{2} \pux^{2\ell-1} \mathcal{H}[f] - \frac{1}{2} \overline{e_0} \pux^{2\ell-1} f \\
& = & (- \overline{e_0} \pux^{2\ell-1}) \left ( \frac{1}{2} f + \frac{1}{2} \overline{e_0} \mathcal{H}[f] \right )  \ = \ (- \overline{e_0} \pux)^{2\ell-1} \left ( \mathcal{A} \mathcal{S} [f] \right ) \\[4mm]
\mathcal{C}_{-2\ell-1}^+ [f]  =  \lim_{x_0 \rightarrow 0+} \mathcal{C}_{-2\ell-1} [f] & = & \frac{1}{2} \pux^{2\ell} f + \frac{1}{2} \overline{e_0} \pux^{2\ell} \mathcal{H} [f] \\
& = & \pux^{2\ell} \left ( \frac{1}{2} f + \frac{1}{2} \overline{e_0} \mathcal{H}[f] \right ) \ = \ \pux^{2\ell} \left ( \mathcal{A} \mathcal{S} [f] \right ) = (- \overline{e_0} \pux)^{2\ell} \left ( \mathcal{A} \mathcal{S} [f] \right )
\end{eqnarray*}
which belong to the Clifford--Hardy space $H^2(\mR^m)$. This
leads to the following commutative schemes for a function $f \in W_2^k(\mR^m)$:
$$
\begin{array}{ccccc}
\mathcal{C}_{-2\ell-1}[f] & \xrightarrow{\hspace*{2mm} - \overline{e_0} \pux \hspace*{2mm}} & \mathcal{C}_{-2\ell-2}[f] & \xrightarrow{\hspace*{2mm} - \overline{e_0} \pux \hspace*{2mm}}  & \mathcal{C}_{-2\ell-3}[f] \\[2mm]
\hspace*{-8mm} ^{\hspace*{0.1mm}_{x_0 \rightarrow 0+}} \downarrow & & \downarrow & & \downarrow \\
(- \overline{e_0} \pux)^{2\ell} \left ( \mathcal{A} \mathcal{S} [f] \right ) & \xrightarrow{\hspace*{2mm} - \overline{e_0} \pux \hspace*{2mm}} & (- \overline{e_0} \pux)^{2\ell+1} \left ( \mathcal{A} \mathcal{S} [f] \right ) & \xrightarrow{\hspace*{2mm} - \overline{e_0} \pux \hspace*{2mm}} & (- \overline{e_0} \pux)^{2\ell+2} \left ( \mathcal{A} \mathcal{S} [f] \right )
\end{array}
$$
The above scheme reflects at the level of the convolution operators, the commutative schemes (\ref{cscheme2}) at the level of the convolution kernels.

\subsection{Explicit expression of the downstream potentials}

In Subsection 4.1 we have already obtained the explicit expressions of the harmonic potentials $A_{-1}, B_{-1}, A_{-2}$ and $B_{-2}$ (see \ref{A-1B-1}, \ref{A-2B-2}). Putting forward for the harmonic potentials $A_{-k}$ and $B_{-k}$ the following form:
$$
A_{-k} = \frac{2}{\sigma_{m+1}} \ \frac{1}{|x|^{m+2k-1}} \ P_k(x_0,|\ux|^2)  \quad {\rm and} \quad 
B_{-k} = \frac{2}{\sigma_{m+1}} \ \frac{\ux}{|x|^{m+2k-1}} \ Q_{k-1}(x_0,|\ux|^2) 
$$
where $P_k$ and $Q_{k-1}$ are scalar--valued homogeneous polynomials of degree $k$ and $k-1$ respectively, it is shown, by a direct calculation, that these polynomials satisfy the following recurrence relations:
$$
P_{k+1}(t,u^{2}) = (t^2+u^{2}) \, \p_t P_k - (m+2k-1) \, t \, P_k  \quad , \quad P_1(t,u^{2}) = t
$$
and
$$
Q_{k}(t,u^{2}) = (t^2+u^{2}) \, \p_t Q_{k-1} - (m+2k-1) \, t \, Q_{k-1}  \quad , \quad Q_0(t,u^{2}) = -1
$$
The fact that $A_{-k}$ and $B_{-k}$ are related by harmonic conjugacy leads to the following intertwined relations for those polynomials:
$$
Q_{k}(x_0,|\ux|^2) = (m+2k-1) \, P_k + \frac{x_0^2 + |\ux|^2}{|\ux|^2} \, \ux \pux P_k
$$
and
$$
P_{k+1}(x_0,|\ux|^2) = \left(m x_0^2 - (2k-1)|\ux|^2\right) Q_{k-1} - (x_0^2 + |\ux|^2) \, \pux Q_{k-1} \, \ux
$$
It is possible to obtain an explicit expression for $P_{k}(t,u^{2})$ and $Q_{k}(t,u^{2})$ in terms of well-known orthogonal polynomials. This is achieved in the following way. First rewrite these polynomials as
\begin{align*}
P_{k}(t,u^{2}) &= u^{k} \widetilde{P}_{k}\left( \frac{t}{u} \right)\\
Q_{k}(t,u^{2}) &= u^{k} \widetilde{Q}_{k}\left( \frac{t}{u} \right)
\end{align*}
with $\widetilde{P}_{k}(w)$ and $\widetilde{Q}_{k}(w)$ polynomials of degree $k$ in the variable $w=t/u$. The recursion relations for $P_{k}$ and $Q_{k}$ may now be rewritten as
\begin{align*}
\widetilde{P}_{k}(w) & = (1+w^{2}) \partial_{w} \widetilde{P}_{k-1}(w) - (m+2k-3) w \widetilde{P}_{k-1}(w), \qquad \widetilde{P}_{0}(w) = -1/(m-1)\\
\widetilde{Q}_{k}(w) & = (1+w^{2}) \partial_{w} \widetilde{Q}_{k-1}(w) - (m+2k-1) w \widetilde{Q}_{k-1}(w), \qquad \widetilde{Q}_{0}(w) =1
\end{align*}
Using the operator identity
\[
(1+ w^{2}) \partial_{w} +2(\alpha+1)w = (1+w^{2})^{-\alpha} \partial_{w} (1+w^{2})^{\alpha+1}
\]
we subsequently find
\begin{align*}
\widetilde{P}_{k}(w) & = - \frac{1}{m-1}(1+w^{2})^{k + \frac{m-1}{2}} (\partial_{w})^{k} (1+w^{2})^{- \frac{m-1}{2}}\\[1mm]
\widetilde{Q}_{k}(w) & = - (1+w^{2})^{k + \frac{m+1}{2}} (\partial_{w})^{k} (1+w^{2})^{- \frac{m+1}{2}}
\end{align*}
Comparing this result with the Rodrigues' formula for the Gegenbauer polynomials, we obtain:
$$
\widetilde{P}_{k}(w) = (-1)^{k+1} 2^k i^k k! \frac{1}{m-1}  \frac{\Gamma(-m-2k+2)}{\Gamma(-m-k+2)} \frac{\Gamma(-\frac{m}{2} + \frac{3}{2})}{\Gamma(-\frac{m}{2} + \frac{3}{2}-k)} C^{-k+1-\frac{m}{2}}_{k}(i w)
$$
and
$$
\widetilde{Q}_{k}(w) = (-1)^{k+1} 2^k i^k k! \frac{\Gamma(-m-2k)}{\Gamma(-m-k)} \frac{\Gamma(-\frac{m}{2} + \frac{1}{2})}{\Gamma(-\frac{m}{2} + \frac{1}{2}-k)} C^{-k-\frac{m}{2}}_{k}(i w)
$$
where $i$ is the imaginary unit, eventually leading to
$$
A_{-k} = (-1)^{k+1} \frac{2^{k+1}}{\sigma_{m+1}} \  k! \ \frac{1}{m-1}  \frac{\Gamma(-m-2k+2)}{\Gamma(-m-k+2)} \frac{\Gamma(-\frac{m}{2} + \frac{3}{2})}{\Gamma(-\frac{m}{2} + \frac{3}{2}-k)} \ \frac{|\ux|^k}{|x|^{m+2k-1}} \ i^k C^{-k+1-\frac{m}{2}}_{k}\left( i \frac{ x_0}{|\ux|} \right)
$$
$$
 = - \frac{2^{2k+1}}{\sigma_{m+1}} \ \frac{1}{m-1}  \frac{\Gamma(-m-2k+2)}{\Gamma(-m-k+2)} \frac{\Gamma(-\frac{m}{2} + \frac{3}{2})}{\Gamma(-\frac{m}{2} + \frac{3}{2}-k)} \frac{\Gamma(-\frac{m}{2}+1)}{\Gamma(-\frac{m}{2}-k+1)}  \ \frac{x_0^k}{|x|^{m+2k-1}} \ _2 F_1\left(-\frac{k}{2}, \frac{1-k}{2}; \frac{m}{2}; - \frac{|\ux|^2}{x_0^2}\right)
 $$
and
$$
B_{-k} = (-1)^{k} \frac{2^{k}}{\sigma_{m+1}} \  (k-1)! \  \frac{\Gamma(-m-2k+2)}{\Gamma(-m-k+1)} \frac{\Gamma(-\frac{m}{2} + \frac{1}{2})}{\Gamma(-\frac{m}{2} + \frac{3}{2}-k)} \ \frac{|\ux|^{k-1} \ux}{|x|^{m+2k-1}} \ i^{k-1} C^{-k+1-\frac{m}{2}}_{k-1}\left( i \frac{ x_0}{|\ux|} \right)
$$
$$
 = - \frac{2^{2k+1}}{\sigma_{m+1}} \  \frac{\Gamma(-m-2k+2)}{\Gamma(-m-k+1)} \frac{\Gamma(-\frac{m}{2} + \frac{1}{2})}{\Gamma(-\frac{m}{2} + \frac{3}{2}-k)} \frac{\Gamma(-\frac{m}{2}+1)}{\Gamma(-\frac{m}{2}-k+2)}  \ \frac{x_0^{k-1} \ux}{|x|^{m+2k-1}} \ _2 F_1\left(-\frac{k+1}{2}, 1-\frac{k}{2}; \frac{m}{2}; - \frac{|\ux|^2}{x_0^2}\right)
 $$
 \\[2mm]
Finally, to give an idea, let us state the explicit expressions of the potentials $A_{-k}$ and $B_{-k}$ for a couple of low values of $k$:
\begin{eqnarray*}
A_{-3} (x_0,|\ux|) & = & \frac{2}{\sigma_{m+1}} \ \frac{1}{|x|^{m+5}}  \left(   m(m+1) x_0^3  - 3(m+1) x_0  |\ux|^2   \right)\\
A_{-4} (x_0,|\ux|) & = & \frac{2}{\sigma_{m+1}} \ \frac{1}{|x|^{m+7}}  \left( - m(m+1)(m+2) x_0^4  + 6(m+1)(m+2) x_0^2  |\ux|^2   - 3(m+1) |\ux|^4  \right)
\end{eqnarray*}
and
\begin{eqnarray*}
B_{-3} (x_0,|\ux|) & = & \frac{2}{\sigma_{m+1}} \ \frac{\ux}{|x|^{m+5}}  \ (m+1) \left(  - (m+2) x_0^2 +  |\ux|^2    \right)\\
B_{-4} (x_0,|\ux|) & = & \frac{2}{\sigma_{m+1}} \ \frac{\ux}{|x|^{m+7}}  \ (m+1)(m+3) \left( (m+2) x_0^3  - 3 x_0  |\ux|^2  \right)
\end{eqnarray*}


\section{Upstream potentials}
\label{upstream}

\subsection{The monogenic logarithmic function}

Recall that Green's function $A_0(x_0,\ux) = - \frac{2}{m-1} \frac{1}{\sigma_{m+1}} \frac{1}{|x|^{m-1}}$, (\ref{A0}), and its conjugate harmonic $B_0(x_0,\ux) =  \frac{2}{\sigma_{m+1}} \frac{\ux}{|\ux|^m} F_m ( \frac{|\ux|}{x_0} )$, (\ref{B0}), satisfy in $\mR^{m+1}_+$ the system
$$
\left \{ \begin{array}{rcl}
\p_{x_0} A_0 & = & - \pux B_0 \ = \ P \ = \ A_{-1} \\[2mm]
\p_{x_0} B_0 & = & - \pux A_0 \ = \ Q \ = \ B_{-1} \\[2mm]
\end{array} \right .
$$
(see also (\ref{CR12}) in Section 3) from which it follows that
\begin{equation}
\overline{D} A_0 = \frac{1}{2} \left ( \p_{x_0} - \overline{e_0} \pux \right ) A_0 \ = \ \frac{1}{2} P + \frac{1}{2} \overline{e_0} Q = C_{-1}
\label{rel1}
\end{equation}
and
\begin{equation}
\overline{D} \overline{e_0} B_0 = \frac{1}{2} \left ( \overline{e_0} \p_{x_0} - \pux \right ) B_0 \ = \ \frac{1}{2} \overline{e_0} Q  + \frac{1}{2} P = C_{-1}
\label{rel2}
\end{equation}
Relations (\ref{rel1}) and (\ref{rel2}) express the fact that $A_0(x_0,\ux)$ and $\overline{e_0} B_0(x_0,\ux)$ are conjugate harmonic potentials (or primitives), with respect to the operator $\overline{D}$, of the Cauchy kernel $C_{-1}(x_0,\ux)$ in $\mR^{m+1}_+$. Putting, as in Section 3, $C_0(x_0,\ux) = \frac{1}{2} A_0(x_0,\ux) + \frac{1}{2} \overline{e_0} B_0 (x_0,\ux)$, it is readily seen that
$$
\overline{D} C_0 (x_0,\ux) = \frac{1}{2} C_{-1} (x_0,\ux) + \frac{1}{2} C_{-1} (x_0,\ux) = C_{-1}(x_0,\ux)
$$
which implies that $C_0(x_0,\ux)$ is a monogenic potential (or primitive), with respect to te operator $\overline{D}$, of the Cauchy kernel $C_{-1}(x_0,\ux)$ in $\mR^{m+1}_+$. Moreover there holds, in view of $DC_0 = \frac{1}{2} (\p_{x_0} +\overline{e_0} \pux) C_0 = 0$, that
$$
C_{-1} = \overline{D} C_0 = \p_{x_0} C_0 = \left ( - \overline{e_0} \pux \right ) C_0
$$
Recall the distributional limits for $x_0 \rightarrow 0+$ of $A_0(x_0,\ux)$ and $B_0(x_0,\ux)$:
\begin{eqnarray*}
a_0(\ux) & = & - \frac{2}{m-1} \frac{1}{\sigma_{m+1}} \mbox{Fp} \frac{1}{|\ux|^{m-1}} \ = \ - \frac{2}{m-1} \frac{1}{\sigma_{m+1}} T^\ast_{-m+1} \\
b_0(\ux) & = & \frac{1}{\sigma_m} \frac{\ux}{|\ux|^m} \ = \ \frac{1}{\pi} \frac{1}{\sigma_m} U^\ast_{-m+1}
\end{eqnarray*} 
(see also Section 3, (\ref{a0}) and (\ref{b0})), yielding
$$
c_0(\ux) = \lim_{x_0 \rightarrow 0+} C_0(x_0,\ux) = - \frac{1}{m-1} \frac{1}{\sigma_{m+1}} \mbox{Fp} \frac{1}{|\ux|^{m-1}} + \frac{1}{2} \overline{e_0} \frac{1}{\sigma_m} \frac{\ux}{|\ux|^m}
$$
As was expected these distributional boundary values are intimately related, as is shown in the following lemma.
\begin{lemma}
\label{lemintiem}
One has
\begin{itemize}
\item[(i)] $-\pux a_0 = b_{-1} = H$
\item[(ii)] $-\pux b_0 = a_{-1} = \delta$
\item[(iii)] $\mathcal{H} \left [a_0 \right ] = b_0$
\item[(iv)] $\mathcal{H} \left [b_0 \right ] = a_0$
\item[(v)] $c_{-1} \ast a_0 = c_{-1} \ast \overline{e_0} b_0 = c_{-1} \ast c_0 = c_0$ 
\item[(vi)] $-\overline{e_0} \pux c_0 = c_{-1}$
\item[(vii)] $ \overline{e_0} \mathcal{H} \left [c_0 \right ] = c_0$ 
\end{itemize}
\end{lemma}

\pf
We make use of the calculation rules for the $T^\ast$-- and $U^\ast$--distributions, recalled in Proposition \ref{prop1}. For (i) we have
$$
- \pux a_0 = \frac{2}{m-1} \frac{1}{\sigma_{m+1}} \pux T^\ast_{-m+1} = \frac{2}{m-1} \frac{1}{\sigma_{m+1}} (-m+1) U^\ast_{-m} = -\frac{2}{\sigma_{m+1}} U^\ast_{-m} = b_{-1} = H
$$
while for (ii) 
$$
- \pux b_0 = -\frac{1}{\pi} \frac{1}{\sigma_m} \pux U^\ast_{-m+1} = -\frac{1}{\pi} \frac{1}{\sigma_m} (-2\pi) T^\ast_{-m} = -\frac{2}{\sigma_m} T^\ast_{-m} = a_{-1} = \delta
$$
Then (iii) is obtained by 
\begin{eqnarray*}
\mathcal{H} \left [ a_0 \right ] &=& H \ast a_0 = \left ( - \frac{2}{\sigma_{m+1}} U^\ast_{-m} \right ) \ast  \left ( - \frac{2}{m-1} \frac{1}{\sigma_{m+1}} T^\ast_{-m+1} \right ) \\
& = & \frac{4}{m-1} \frac{1}{\left ( \sigma_{m+1} \right )^2} U^\ast_{-m} \ast T^\ast_{-m+1} \ = \ \frac{4}{m-1} \frac{1}{
\left ( \sigma_{m+1} \right )^2} \pi^\frac{m}{2} \frac{\Gamma \left ( \frac{m}{2} \right )}{\Gamma \left ( \frac{m+1}{2} \right ) \Gamma \left ( \frac{m-1}{2}\right )} U^\ast_{-m+1} \\
& = & \frac{\Gamma \left ( \frac{m}{2} \right ) }{2 \pi^{\frac{m}{2}+1} } U^\ast_{-m+1} \ = \ \frac{1}{\pi} \frac{1}{\sigma_m} U^\ast_{-m+1} \ = \ b_0
\end{eqnarray*}
from which also (iv) follows: $\mathcal{H} \left [ b_0 \right ] = \mathcal{H}^2 \left [ a_0 \right ] = a_0$. 
To obtain (v) it suffices to observe that
\begin{eqnarray*}
c_{-1} \ast a_0 & = & \left ( \frac{1}{2} a_{-1} + \frac{1}{2} \overline{e_0} b_{-1} \right ) \ast a_0 \ = \ \left ( \frac{1}{2} \delta + \frac{1}{2} \overline{e_0} H \right ) \ast a_0 \ = \ \frac{1}{2} a_0 + \frac{1}{2} \overline{e_0} b_0 \ = \  c_0 \\
c_{-1} \ast \overline{e_0} b_0 & = & \left ( \frac{1}{2} \delta + \frac{1}{2} \overline{e_0} H \right ) \ast \overline{e_0} b_0 \ = \ \frac{1}{2} \overline{e_0} b_0 + \frac{1}{2} \mathcal{H} \left [ b_0 \right ] \ = \ \frac{1}{2} \overline{e_0} b_0 + \frac{1}{2} a_0 \ = \ c_0 \\
c_{-1} \ast c_0 & = & c_{-1} \ast \left ( \frac{1}{2} a_0 + \frac{1}{2} \overline{e_0} b_0 \right ) \ = \ \frac{1}{2} c_0 + \frac{1}{2} c_0 \ = \ c_0
\end{eqnarray*}
while (vi) is directly obtained by
$$
- \overline{e_0} \pux c_0 = \overline{e_0} \left ( -\pux \frac{1}{2} a_0 - \pux \frac{1}{2} \overline{e_0} b_0 \right ) \ = \ \frac{1}{2} \overline{e_0} b_{-1} + \frac{1}{2} a_{-1} = c_{-1}
$$
Finally, we have
$$
\mathcal{H} \left [ c_0 \right ] = \frac{1}{2} \mathcal{H} \left [ a_0 \right ] + \frac{1}{2} \mathcal{H} \left [ \overline{e_0} b_0 \right ] = \frac{1}{2} b_0 + \frac{1}{2} e_0 a_0
$$
from which (vii) follows:
$$
\overline{e_0} \mathcal{H} \left [ c_0 \right ] = \frac{1}{2} a_0 + \frac{1}{2} \overline{e_0} b_0 = c_0
$$
\qed
Note the following commutative schemes which are each others Hilbert image:
\begin{equation}
\begin{array}{ccc}
A_0(x_0,\ux) & \xrightarrow{\hspace*{1mm} - \pux \hspace*{1mm}} & B_{-1}(x_0,\ux) = Q(x_0,\ux) \\[2mm]
\hspace*{-8mm} ^{\hspace*{0.1mm}_{x_0 \rightarrow 0+}} \downarrow & & \downarrow \\
a_{0}(\ux) &  \xrightarrow{\hspace*{1mm} - \pux \hspace*{1mm}} & b_{-1}(\ux) = H(\ux)
\end{array}
\quad \mbox{and} \qquad 
\begin{array}{ccc}
B_{0}(x_0,\ux) &  \xrightarrow{\hspace*{1mm} - \pux \hspace*{1mm}} & A_{-1}(x_0,\ux) = P(x_0,\ux) \\[2mm]
\hspace*{-7mm} ^{\hspace*{0.1mm}_{x_0 \rightarrow 0+}} \downarrow & & \downarrow \\
b_{0}(\ux) &  \xrightarrow{\hspace*{1mm} - \pux \hspace*{1mm}} & a_{-1}(\ux) = \delta(\ux)
\end{array}
\label{cschemes4}
\end{equation}
The following lemma, also in terms of distributions, makes the relationships between the harmonic potentials $A_0$, $B_0$, $A_{-1} = P$ and $B_{-1}=Q$ more transparent.
\begin{lemma}
\label{lem52}
In distributional sense one has, convolutions being taken in the variable $\ux \in \mR^m$:
\begin{itemize}
\item[(i)] $A_0 = a_0 \ast A_{-1} = b_0 \ast B_{-1}$
\item[(ii)] $B_0 = b_0 \ast A_{-1} = a_0 \ast B_{-1}$
\item[(iii)] $H \ast A_0 = B_0 = A_0 \ast H$
\item[(iv)] $H \ast B_0 = A_0 = B_0 \ast H$
\end{itemize}
\end{lemma}

\pf
\noindent (i)(ii) The technique is the same as the one used in the proof of Proposition \ref{propdownstream}.\\[-2mm]

\noindent (iii)(iv) It suffices to observe that 
$$
A_0 \ast H = A_{-1} \ast a_0 \ast H = A_{-1} \ast b_0 = B_0
$$
and
$$
B_0 \ast H = A_0 \ast H \ast H = A_0
$$
\qed

Similar properties hold for the monogenic potentials $C_0$ and $C_{-1}$.
\begin{lemma}
\label{lem53}
In distributional sense one has, convolutions being taken in the variable $\ux \in \mR^m$:
\begin{itemize}
\item[(i)] $C_0 = C_{-1} \ast a_0 = C_{-1} \ast \overline{e_0} b_0 = C_{-1} \ast c_0$
\item[(ii)] $C_0 = A_{-1} \ast c_0 = \overline{e_0} B_{-1} \ast c_0$
\end{itemize}
\end{lemma}

\pf
Making use of the results of Lemma \ref{lem52}, we have for (i)
\begin{eqnarray*}
C_{-1} \ast a_0 & = & \left ( \frac{1}{2} A_{-1} + \frac{1}{2} \overline{e_0} B_{-1} \right ) \ast a_0 \ = \ \frac{1}{2} A_0 + \frac{1}{2} \overline{e_0} B_0 \ = \ C_0 \\
C_{-1} \ast \overline{e_0} b_0 & = & \left ( \frac{1}{2} A_{-1} + \frac{1}{2} \overline{e_0} B_{-1} \right ) \ast \overline{e_0} b_0 \ = \ \frac{1}{2} \overline{e_0} B_0 + \frac{1}{2} A_0 \ = \ C_0 \\
C_{-1} \ast c_0 & = & C_{-1} \ast \left ( \frac{1}{2} a_0 + \frac{1}{2} \overline{e_0} b_0 \right ) = \frac{1}{2} C_0 + \frac{1}{2} C_0 \ = \ C_0
\end{eqnarray*}
while for (ii)
\begin{eqnarray*}
A_{-1} \ast c_0 & = & A_{-1} \ast \left ( \frac{1}{2} a_0 + \frac{1}{2} \overline{e_0} b_0 \right ) \ = \ \frac{1}{2} A_0 + \frac{1}{2} \overline{e_0} B_0 \ = \ C_0 \\
\overline{e_0} B_{-1} \ast c_0 & = & \overline{e_0} B_{-1} \ast \left ( \frac{1}{2} a_0 + \frac{1}{2} \overline{e_0} b_0 \right ) \ = \ \frac{1}{2} \overline{e_0} B_0 + \frac{1}{2} B_{-1} \ast b_0 \ = \ C_0
\end{eqnarray*}
\qed

\begin{remark}
The results (i)--(ii) of Lemma \ref{lem52} are the analogues of the results of Proposition \ref{propdownstream} for the case where $j=-1$, $k=1$.
\end{remark}

The potential kernels $A_0(x_0,\ux)$, $B_0(x_0,\ux)$ and $C_0(x_0,\ux)$ may now be used in their corresponding convolution operators defining the conjugate harmonic potentials in $\mR^{m+1}_+$
$$
\mathcal{A}_0[f](x_0,\ux) = A_0(x_0,\cdot) \ast f(\cdot)(\ux) \quad \mbox{and} \quad \mathcal{B}_0[f](x_0,\ux) = B_0(x_0,\cdot) \ast f(\cdot)(\ux)
$$
and the monogenic potential
$$
\mathcal{C}_0[f](x_0,\ux) = C_0(x_0,\cdot) \ast f(\cdot) (\ux) = \frac{1}{2} \mathcal{A}_0[f](x_0,\ux) + \frac{1}{2} \overline{e_0} \mathcal{B}_0[f](x_0,\ux)
$$
The properties they enjoy, summarized in the next proposition, reflect the corresponding properties of the potential kernels.
\begin{proposition}
\label{prop51}
For a Schwartz function or a distribution $f$ one has
\begin{itemize}
\item[(i)] $\overline{D} \mathcal{A}_0[f] = \overline{D} \mathcal{B}_0[f] = \overline{D} \mathcal{C}_0[f] = C_{-1}[f]$
\item[(ii)] $\lim_{x_0 \rightarrow 0+} \mathcal{A}_0[f] = a_0 \ast f = - (- \Delta)^{-\frac{1}{2}} [f]$\\
$\lim_{x_0 \rightarrow 0+} \mathcal{B}_0[f] = b_0 \ast f = - E \ast f = T[f]$ \\
$\lim_{x_0 \rightarrow 0+} \mathcal{C}_0[f] = c_0 \ast f = a_0 \ast \mathcal{A} \mathcal{S}[f] = \overline{e_0} b_0 \ast f$
\item[(iii)] $\overline{e_0} \mathcal{H} \left [c_0 \ast f \right ] = c_0 \ast f$
\item[(iv)] $\mathcal{C}_0[f] = \mathcal{C}_{-1} \left [a_0 \ast f \right ] = \mathcal{C}_{-1} \left [ b_0 \ast \mathcal{H}[f] \right ] = \mathcal{C}_{-1} \left [\overline{e_0} b_0 \ast f \right ] = \mathcal{C}_{-1} \left [ a_0 \ast \overline{e_0} \mathcal{H}[f] \right ]$ \\
$\phantom{\mathcal{C}_0[f]} = \mathcal{C}_{-1} \left [ c_0 \ast f \right ] = \mathcal{P} \left [ c_0 \ast f \right ] = \mathcal{A}_0 \left [ \mathcal{A}\mathcal{S}[f] \right ] = \overline{e_0} \mathcal{B}_0 \left [ \mathcal{A} \mathcal{S}[f] \right ]$
\end{itemize}
\end{proposition}
Note also the following commutative scheme, which may be derived from the commutative schemes (\ref{cschemes4}):
$$
\begin{array}{rcl}
\mathcal{C}_{0}[f] & \xrightarrow{\hspace*{2mm} -\overline{e_0} \pux \hspace*{2mm}} & \mathcal{C}_{-1}[f] \\[2mm]
^{\hspace*{0.1mm}_{x_0 \rightarrow 0+}} \downarrow \hspace*{4mm} & & \hspace*{4mm} \downarrow ^{\hspace*{0.1mm}_{x_0 \rightarrow 0+}}\\
c_0 \ast f & \xrightarrow{\hspace*{2mm} -\overline{e_0} \pux \hspace*{2mm}} & c_{-1} \ast f =  \mathcal{A} \mathcal{S}[f] 
\end{array}
$$

\begin{remark}
As already explained in the introduction, in the upper half of the complex plane the function $\ln(z)$ is a holomorphic potential (or primitive) of the Cauchy kernel $\frac{1}{z}$ and its real and imaginary components are the fundamental solution $\ln |z|$ of the Laplace operator, and its conjugate harmonic $i\, {\rm arg}(z)$ respectively. By similarity we could say that $C_0(x_0,\ux) = \frac{1}{2} A_0(x_0,\ux) + \frac{1}{2} \overline{e_0} B_0(x_0,\ux)$, being a monogenic potential of the Cauchy kernel $C_{-1}(x_0,\ux)$ and the sum of the fundamental solution $A_0(x_0,\ux)$ of the Laplace operator and its conjugate harmonic $\overline{e_0} B_0(x_0,\ux)$, is a {\em monogenic logarithmic function} in the upper half--space $\mR^{m+1}_+$.
\end{remark}

\subsection{The potentials of the logarithmic monogenic function}

\noindent
Inspired by the properties contained in Proposition \ref{propdownstream} and Lemma \ref{lem52} we proceed as follows for the construction of harmonic and monogenic potentials of $C_0(x_0,\ux)$ in $\mR^{m+1}_+$. We put
$$
\left \{ \begin{array}{rclcl}
A_1(x_0,\ux) & = & a_0(\cdot) \ast A_0(x_0,\cdot)(\ux) & = & b_0(\cdot) \ast B_0(x_0,\cdot) \\[2mm]
B_1(x_0,\ux) & = & a_0(\cdot) \ast B_0(x_0,\cdot)(\ux) & = & b_0(\cdot) \ast A_0(x_0,\cdot)
\end{array} \right .
$$
and we verify at once that
\begin{equation}
\begin{array}{rcl}
\p_{x_0} A_1 & = & a_0 \ast \p_{x_0} A_0 \ = \ a_0 \ast A_{-1} \ = \ A_{-1} \ast a_0 \ = \ \mathcal{P}[a_0] \ = \ A_0\\[2mm]
\p_{x_0} A_1 & = & b_0 \ast \p_{x_0} B_0 \ = \ b_0 \ast (-\pux A_{0}) \ = \ (-b_0 \pux) \ast A_0 \ = \ \delta \ast A_0 \ = \ A_0
\end{array}
\label{verify1}
\end{equation}
and
\begin{equation}
\begin{array}{rcl}
\p_{x_0} B_1 & = & a_0 \ast \p_{x_0} B_0 \ = \ (-a_0 \pux) \ast A_{0} \ = \ \mathcal{H}[A_0] \ = \ B_0\\[2mm]
\p_{x_0} B_1 & = & b_0 \ast \p_{x_0} A_0 \ = \ b_0 \ast (-\pux B_{0}) \ = \ (-b_0 \pux) \ast B_0 \ = \ \delta \ast B_0 \ = \ B_0
\end{array}
\label{verify2}
\end{equation}
while also
\begin{equation}
\begin{array}{rcl}
- \pux A_1 & = & - \pux a_0 \ast A_0 \ = \ H \ast A_0 \ = \ \mathcal{H}[A_0] \ = \ B_0 \\[2mm]
- \pux A_1 & = & - \pux b_0 \ast B_0 \ = \ \delta \ast B_0 \ = \ B_0 
\end{array}
\label{verify3}
\end{equation}
and
\begin{equation}
\begin{array}{rcl}
- \pux B_1 & = & - \pux a_0 \ast B_0 \ = \ \mathcal{H}[B_0] \ = \ A_0 \\[2mm]
- \pux B_1 & = & - \pux b_0 \ast A_0 \ = \ \delta \ast A_0 \ = \ A_0 
\end{array}
\label{verify4}
\end{equation}
The relations (\ref{verify1})--(\ref{verify4}) precisely are the relations needed for $A_1(x_0,\ux)$ and $B_1(x_0,\ux)$ to be conjugate harmonic potentials in $\mR^{m+1}_+$ of the function $C_0(x_0,\ux)$. It then follows at once that
$$
C_1(x_0,\ux) = \frac{1}{2} A_1(x_0,\ux) + \frac{1}{2} \overline{e_0} B_1(x_0,\ux)
$$
is a monogenic potential in $\mR^{m+1}_+$ of $C_0$ and there holds that $\overline{D} C_1 = \p_{x_0} C_1 = (-\overline{e_0} \pux) C_1 = C_0$. Also note that the conjugate harmonic potentials $A_1(x_0,\ux)$ and $B_1(x_0,\ux)$ form a Hilbert pair, the Hilbert transform being taken in the variable $\ux \in \mR^m$. Indeed, we have
\begin{eqnarray*}
\mathcal{H} \left [ A_1 \right ] &=& H(\cdot) \ast A_1(x_0,\cdot) (\ux) \ = \ H(\cdot) \ast a_0(\cdot) \ast A_0(x_0,\cdot)(\ux) \ = \ b_0(\cdot) \ast A_0(x_0,\cdot) \ = \ B_1(x_0,\ux) \\
\mathcal{H} \left [ B_1 \right ] & = & \mathcal{H}^2 \left [ A_1 \right ] \ = \ A_1
\end{eqnarray*}
The distributional limits for $x_0 \rightarrow 0+$ of the conjugate harmonic potentials $A_1$ and $B_1$ are given by
$$
\left \{ \begin{array}{rcl}
a_1(\ux) & = & \lim_{x_0 \rightarrow 0+} A_1(x_0,\ux) \ = \ a_0(\cdot) \ast a_0(\cdot)(\ux) \ = \ b_0(\cdot) \ast b_0(\cdot) (\ux) \\[2mm]
b_1(\ux) & = & \lim_{x_0 \rightarrow 0+} B_1(x_0,\ux) \ = \ a_0(\cdot) \ast b_0(\cdot)(\ux) \ = \ b_0(\cdot) \ast a_0(\cdot) (\ux)
\end{array} \right .
$$
Making use of the calculation rules for the convolution of the $T^\ast$-- and $U^\ast$--distributions (see Section 2, Proposition 2.1), these distributional boundary values are explicitly given by
\begin{equation}
\label{a1b1}
\left \{ \begin{array}{rcl}
a_1(\ux) & = & \phantom{-} \displaystyle\frac{1}{\pi} \displaystyle\frac{1}{\sigma_m} \displaystyle\frac{1}{m-2} \, T^\ast_{-m+2} \ = \ \displaystyle\frac{1}{\sigma_m} \displaystyle\frac{1}{m-2} \displaystyle\frac{1}{|\ux|^{m-2}}
 \\[5mm]
b_1(\ux) & = & - \displaystyle\frac{1}{\pi} \displaystyle\frac{1}{\sigma_{m+1}} \displaystyle\frac{1}{m-1} \, U^\ast_{-m+2} \ = \ - \displaystyle\frac{1}{\sigma_{m+1}} \displaystyle\frac{2}{m-1} \displaystyle\frac{\ux}{|\ux|^{m-1}}
\end{array} \right .
\end{equation}

They show the following properties, where we have put, quite naturally, $c_1(\ux) = \frac{1}{2} a_1(\ux) + \frac{1}{2} \overline{e_0} b_1(\ux)$.
\begin{lemma}
\label{lem54}
\rule{0mm}{0mm}
\begin{itemize}
\item[(i)] $- \pux a_1 = b_0$, $-\pux b_1 = a_0$, $- \overline{e_0} \pux c_1 = c_0$
\item[(ii)] $\mathcal{H} \left [ a_1 \right ] = b_1$, $\mathcal{H} \left [ b_1 \right ] = a_1$, $\overline{e_0} \mathcal{H} \left [ c_1 \right ] = c_1$
\item[(iii)] $c_{-1} \ast a_1 = c_{-1} \ast \overline{e_0} b_1 = c_{-1} \ast c_1 = c_1$
\item[(iv)] $a_0 \ast c_0 = c_0 \ast a_0 = c_1$, $\overline{e_0} b_0 \ast c_0 = c_0 \ast \overline{e_0} b_0 = c_1$
\end{itemize}
\end{lemma}

\pf
For (i) we consecutively have
\begin{eqnarray*}
- \pux a_1 & = & - \frac{1}{\pi} \frac{1}{\sigma_m} \frac{1}{m-2} \, \pux T^\ast_{-m+2} =  - \frac{1}{\pi} \frac{1}{\sigma_m} \frac{1}{m-2} (-m+2) \, U^\ast_{-m+1} = \frac{1}{\pi} \frac{1}{\sigma_m} \, U^\ast_{-m+1} = b_0 \\
- \pux b_1 & = & \frac{1}{\pi} \frac{1}{\sigma_{m+1}} \frac{1}{m-1} \, \pux U^\ast_{-m+2} =  \frac{1}{\pi} \frac{1}{\sigma_{m+1}} \frac{1}{m-1} (-2\pi) \, T^\ast_{-m+1} = - \frac{2}{m-1} \frac{1}{\sigma_{m+1}} \, T^\ast_{-m+1} = a_0 \\
\end{eqnarray*}
and
$$
- \overline{e_0} \pux c_1 = - \overline{e_0} \pux \left ( \frac{1}{2} a_1 + \frac{1}{2} \overline{e_0} b_1 \right ) = \overline{e_0} \frac{1}{2} b_0 + \frac{1}{2} a_0 = c_0
$$
For (ii) we obtain
\begin{eqnarray*}
\mathcal{H} \left [ a_1 \right ] & = & \mathcal{H} \left [ a_0 \ast a_0 \right ] = b_0 \ast a_0 = b_1 \\
\mathcal{H} \left [ b_1 \right ] & = & \mathcal{H}^2 \left [ a_1 \right ] = a_1
\end{eqnarray*}
and 
$$
\overline{e_0} \mathcal{H} \left [ c_1 \right ] = \overline{e_0} \left ( \frac{1}{2} \mathcal{H} \left [ a_1 \right ] + \frac{1}{2} \mathcal{H} \left [ \overline{e_0} b_1 \right ] \right ) = \frac{1}{2} \overline{e_0} b_1 + \frac{1}{2} a_1 = c_1
$$
For (iii) it holds that
\begin{eqnarray*}
c_{-1} \ast a_1 & = & \left ( \frac{1}{2} \delta + \frac{1}{2} \overline{e_0} H \right ) \ast a_1 = \frac{1}{2} a_1 + \frac{1}{2} \overline{e_0} b_1 = c_1 \\
c_{-1} \ast \overline{e_0} b_1 & = & \left ( \frac{1}{2} \delta + \frac{1}{2} \overline{e_0} H \right ) \ast (\overline{e_0} b_1 ) = \frac{1}{2} \overline{e_0} b_1 + \frac{1}{2} a_1 = c_1
\end{eqnarray*}
and
$$
c_{-1} \ast c_1 = c_{-1} \ast \left ( \frac{1}{2} a_1 + \frac{1}{2} \overline{e_0} b_1 \right ) = \frac{1}{2} c_1 + \frac{1}{2} c_1 = c_1
$$
Finally, (iv) follows from the following calculations:
\begin{eqnarray*}
a_0 \ast c_0 & = & a_0 \ast \left ( \frac{1}{2} a_0 + \frac{1}{2} \overline{e_0} b_0 \right ) = \frac{1}{2} a_1 + \frac{1}{2} \overline{e_0} b_1 = c_1 \\
c_0 \ast a_0 & = & \left ( \frac{1}{2} a_0 + \frac{1}{2} \overline{e_0} b_0 \right ) \ast a_0 = \frac{1}{2} a_1  + \frac{1}{2} \overline{e_0} b_1 = c_1 
\end{eqnarray*}
and
\begin{eqnarray*}
\overline{e_0} b_0 \ast c_0 & = & \overline{e_0} b_0 \ast \left ( \frac{1}{2} a_0 + \frac{1}{2} \overline{e_0} b_0 \right ) = \frac{1}{2} \overline{e_0} b_1 + \frac{1}{2} a_1 = c_1 \\
c_0 \ast \overline{e_0} b_0 & = & \left ( \frac{1}{2} a_0 + \frac{1}{2} \overline{e_0} b_0 \right ) \ast \overline{e_0} b_0 = \frac{1}{2} \overline{e_0} b_1 + \frac{1}{2} a_1 = c_1
\end{eqnarray*}
\qed

Now that we have the distributional boundary values $a_1(\ux)$ and $b_1(\ux)$ at our disposal, the following relations between the harmonic potentials $(A_1,B_1)$ and $(A_{-1},B_{-1})$ may be readily shown.
\begin{lemma}
\label{lem55}
One has, convolutions being taken in the variable $\ux \in \mR^m$:
\begin{itemize}
\item[(i)] $A_1(x_0,\ux) = a_1(\cdot) \ast A_{-1}(x_0,\cdot)(\ux) = b_1(\cdot) \ast B_{-1} (x_0,\cdot)(\ux)$
\item[(ii)] $B_1(x_0,\ux) = a_1(\cdot) \ast B_{-1}(x_0,\cdot)(\ux) = b_1(\cdot) \ast A_{-1} (x_0,\cdot)(\ux)$
\end{itemize}
\end{lemma}

Note that $( - \Delta_m) a_1 = ( - \pux)^2 a_1 = ( - \pux) b_0 = \delta$, and indeed, in $-a_1 = - \frac{1}{\sigma_m} \frac{1}{m-2} \frac{1}{|\ux|^{m-2}}$ we recognize the fundamental solution of the Laplace operator $\Delta_m$ in $\mR^m$. Also note the commutative schemes which are each others Hilbert image:
\begin{equation}
\begin{array}{ccc}
A_1(x_0,\ux) & \xrightarrow{\hspace*{2mm} - \pux \hspace*{2mm}} & B_{0}(x_0,\ux) \\[2mm]
\hspace*{-8mm} ^{\hspace*{0.1mm}_{x_0 \rightarrow 0+}} \downarrow & & \downarrow \\
a_{1}(\ux) &  \xrightarrow{\hspace*{1mm} - \pux \hspace*{1mm}} & b_{0}(\ux) 
\end{array}
\quad \mbox{and} \qquad 
\begin{array}{ccc}
B_{1}(x_0,\ux) &  \xrightarrow{\hspace*{1mm} - \pux \hspace*{1mm}} & A_{0}(x_0,\ux)  \\[2mm]
\hspace*{-7mm} ^{\hspace*{0.1mm}_{x_0 \rightarrow 0+}} \downarrow & & \downarrow \\
b_{1}(\ux) &  \xrightarrow{\hspace*{1mm} - \pux \hspace*{1mm}} & a_{0}(\ux) 
\end{array}
\label{cschemes6}
\end{equation}
As before, the potential kernels $A_1(x_0,\ux)$, $B_1(x_0,\ux)$ and $C_1(x_0,\ux)$ may be used as convolution kernels to define conjugate harmonic functions and their monogenic sum in $\mR^{m+1}_+$, by putting, for a function $f \in L_2(\mR^m)$:
$$
\left \{ \begin{array}{rcl}
\mathcal{A}_1[f] & = & A_1(x_0,\cdot) \ast f(\cdot) (\ux) \ = \ \mathcal{A}_0 \left [ a_0 \ast f \right ] \ = \ \mathcal{B}_0 \left [ b_0 \ast f \right ] \\[2mm]
\mathcal{B}_1[f] & = & B_1(x_0,\cdot) \ast f(\cdot) (\ux) \ = \ \mathcal{A}_0 \left [ b_0 \ast f \right ] \ = \ \mathcal{B}_0 \left [ a_0 \ast f \right ] 
\end{array} \right .
$$
and
$$
\mathcal{C}_1[f] = C_1(x_0,\cdot) \ast f(\cdot) (\ux) = \mathcal{C}_0 \left [ a_0 \ast f \right ] = \mathcal{C}_0 \left [ \overline{e_0} b_0 \ast f \right ]
$$
The corresponding non--tangential $L_2$--boundary values for $x_0 \rightarrow 0+$ are given by $\lim_{x_0 \rightarrow 0+} \mathcal{A}_1 [f] = a_1 \ast f$, $\lim_{x_0 \rightarrow 0+} \mathcal{B}_1 [f] = b_1 \ast f$ and $\lim_{x_0 \rightarrow 0+} \mathcal{C}_1 [f] = c_1 \ast f$, and the commutative schemes (\ref{cschemes6}) eventually lead to the following one:
$$
\begin{array}{rcl}
\mathcal{C}_{1}[f] & \xrightarrow{\hspace*{2mm} -\overline{e_0} \pux \hspace*{2mm}} & \mathcal{C}_{0}[f] \\[2mm]
^{\hspace*{0.1mm}_{x_0 \rightarrow 0+}} \downarrow \hspace*{4mm} & & \hspace*{4mm} \downarrow ^{\hspace*{0.1mm}_{x_0 \rightarrow 0+}}\\
c_1 \ast f & \xrightarrow{\hspace*{2mm} -\overline{e_0} \pux \hspace*{2mm}} & c_{0} \ast f  
\end{array}
$$
The conjugate harmonic potentials $A_1(x_0,\ux)$ and $B_1(x_0,\ux)$ may now be determined explicitly by a computation similar to the one used in Section 3 to determine the conjugate harmonic of Green's function. For this and subsequent calculations the dimension $m$ is assumed to be great enough in order that the expressions obtained should remain valid. Starting from the equation (\ref{verify1})
$$
\p_{x_0} A_1(x_0,\ux) = A_0(x_0,\ux) = - \frac{2}{m-1} \frac{1}{\sigma_{m+1}} \frac{1}{|x_0 e_0 + \ux|^{m-1}}
$$
we find in $\mR^{m+1}_+$
\begin{equation}
A_1(x_0,\ux) = a_1(\ux) - \frac{2}{m-1} \frac{1}{\sigma_{m+1}} \frac{1}{|\ux|^{m-2}} \, \widetilde{F}_{m-2} \left ( \frac{x_0}{|\ux|} \right )
\label{A1}
\end{equation}
where we have put
$$
\widetilde{F}_{m-2} \left ( u \right ) = \int_0^u \frac{d\zeta}{\left ( 1+ \zeta^2 \right )^\frac{m-1}{2}}
$$
A priori it is not clear that $A_1(x_0,\ux)$ is well--defined for $\ux=0$. However, in virtue of the relation
$$
\widetilde{F}_{m-2}(u) = F_{m-2}(+\infty) - F_{m-2} \left ( \frac{1}{u} \right ) = \frac{\sqrt{\pi}}{2} \frac{2}{m-2} \frac{
\Gamma \left ( \frac{m}{2} \right )}{\Gamma \left ( \frac{m-1}{2} \right )} - F_{m-2} \left ( \frac{1}{u} \right )
$$
expression (\ref{A1}) for $A_1(x_0,\ux)$ is turned into, with $m>2$:
$$
A_1(x_0,\ux) = \frac{2}{m-1} \frac{1}{\sigma_{m+1}} \frac{1}{|\ux|^{m-2}} \, F_{m-2} \left ( \frac{|\ux|}{x_0} \right )
$$ 
or, introducing again the hypergeometric function $_2F_1$,
$$
A_1(x_0,\ux) = \frac{2}{m-1} \frac{1}{\sigma_{m+1}} \frac{1}{m-2} \frac{1}{x_0^{m-2}} \; _2F_1 \left ( \frac{m}{2} -1; \frac{m-1}{2}; \frac{m}{2} ; - \frac{|\ux|^2}{x_0^2} \right )
$$
showing that $A_1(x_0,\ux)$ indeed is well--defined for $|\ux|=0$ with
$$
A_1(x_0,0) = \frac{2}{(m-1)(m-2)} \frac{1}{x_0^{m-2}}, \qquad x_0 > 0
$$
By some lengthy calculations it may be verified that the above function $A_1(x_0,\ux)$ also satisfies the equation
$$
- \pux A_1(x_0,\ux) = B_0(x_0,\ux) = \frac{2}{\sigma_{m+1}} \frac{\ux}{|\ux|^m} \, F_m \left ( \frac{|\ux|}{x_0} \right )
$$
and also shows the distributional limit (\ref{a1b1}) given by
$$
\lim_{x_0 \rightarrow 0+} \ A_1(x_0, \ux) = \frac{1}{\sigma_m} \frac{1}{m-2} \frac{1}{|\ux|^{m-2}} = a_1(\ux)
$$
It is perhaps interesting to mention that in the course of these calculations, use has been made of the following recurrence relation for the function $F_m$:
$$
F_m(v) = \frac{m-2}{m-1} \, F_{m-2}(v) - \frac{1}{m-1} \frac{v^{m-2}}{\left ( 1+v^2 \right )^\frac{m-1}{2}}
$$
For the harmonic potential $B_1(x_0,\ux)$ we start the computation  from equation (\ref{verify2}):
$$
\p_{x_0} B_1 (x_0,\ux) = B_0(x_0,\ux) = \frac{2}{\sigma_{m+1}} \frac{\ux}{|\ux|^m} \, F_m \left ( \frac{|\ux|}{x_0} \right )
$$
leading to the expression, with $m>1$:
$$
B_1(x_0,\ux) = \frac{2}{\sigma_{m+1}} \frac{x_0 \ux}{|\ux|^m} \, F_m \left ( \frac{|\ux|}{x_0} \right ) - \frac{2}{\sigma_{m+1}} \frac{1}{m-1} \frac{\ux}{|x|^{m-1}}
$$
or
$$
B_1(x_0,\ux) = \frac{2}{m} \frac{1}{\sigma_{m+1}} \frac{\ux}{x_0^{m-1}} \, _2F_1 \left ( \frac{m}{2}; \frac{m+1}{2} ; \frac{m}{2}+1 ; - \frac{|\ux|^2}{x_0^2} \right ) - \frac{2}{\sigma_{m+1}} \frac{1}{m-1} \frac{\ux}{|x|^{m-1}}
$$
showing that $B_1(x_0,\ux)$ is well--defined for $\ux=0$ with 
$$
B_1(x_0,0) =0 ,\qquad x_0 > 0
$$
Note that the distributional limit $b_1(\ux)$ is indeed recovered:
$$
\lim_{x_0 \rightarrow 0+} B_1(x_0,\ux) = - \frac{2}{\sigma_{m+1}} \frac{1}{m-1} \frac{\ux}{|x|^{m-1}} = b_1((\ux)
$$
It may now readily be verified that the above function $B_1(x_0,\ux)$ also satisfies equation (\ref{verify4}):
$$
- \pux B_1(x_0,\ux) = A_0(x_0,\ux) = - \frac{2}{m-1} \frac{1}{\sigma_{m+1}} \frac{1}{|x|^{m-1}}
$$

\subsection{The next step}

\noindent
Proceeding in a similar way we may define in the next step
$$
\left \{
\begin{array}{rcl}
A_2(x_0,\ux) & = & a_0(\cdot) \ast A_1(x_0,\cdot)(\ux) = b_0(\cdot) \ast B_1(x_0,\cdot)(\ux) \\[2mm]
B_2(x_0,\ux) & = & a_0(\cdot) \ast B_1(x_0,\cdot)(\ux) = b_0(\cdot) \ast A_1(x_0,\cdot)(\ux)
\end{array}
\right .
$$
and verify that in $\mR^{m+1}_+$ it holds that
\begin{eqnarray*}
\p_{x_0} A_2 &=& a_0 \ast \p_{x_0} A_1 \ = \ a_0 \ast A_0 \ = \ A_1 \\
\p_{x_0} A_2 &=& b_0 \ast \p_{x_0} B_1 \ = \ b_0 \ast B_0 \ = \ A_1
\end{eqnarray*}
and
\begin{eqnarray*}
\p_{x_0} B_2 &=& a_0 \ast \p_{x_0} B_1 \ = \ a_0 \ast B_0 \ = \ B_1 \\
\p_{x_0} B_2 &=& b_0 \ast \p_{x_0} A_1 \ = \ b_0 \ast A_0 \ = \ B_1
\end{eqnarray*}
while also
\begin{eqnarray*}
-\pux A_2 &=& -\pux a_0 \ast A_1 \ = \ H \ast A_1 \ = \ \mathcal{H} \left [ A_1 \right ] = B_1\\
-\pux A_2 &=& -\pux b_0 \ast B_1 \ = \ \delta \ast B_1 \ = \ B_1
\end{eqnarray*}
and
\begin{eqnarray*}
-\pux B_2 &=& -\pux a_0 \ast B_1 \ = \ H \ast B_1 \ = \ \mathcal{H} \left [ B_1 \right ] = A_1\\
-\pux B_2 &=& -\pux b_0 \ast A_1 \ = \ \delta \ast A_1 \ = \ A_1
\end{eqnarray*}
These relations justify $A_2(x_0,\ux)$ and $B_2(x_0,\ux)$ to be called conjugate harmonic potentials in $\mR^{m+1}_+$ of the function $C_1(x_0,\ux)$. It follows that
$$
C_2(x_0,\ux) = \frac{1}{2} A_2(x_0,\ux) + \frac{1}{2} \overline{e_0} B_2(x_0,\ux)
$$
is a monogenic potential in $\mR^{m+1}_+$ of $C_1$ and there also holds $\overline{D} C_2 = \p_{x_0} C_2 = (-\overline{e_0} \pux) C_2 = C_1$. As before, the conjugate harmonic potentials $A_2(x_0,\ux)$ and $B_2(x_0,\ux)$ form a Hilbert pair in the variable $\ux \in \mR^m$:
$$
\left \{ \begin{array}{rcl}
\mathcal{H} \left [ A_2(x_0,\ux) \right ] & = & H(\cdot) \ast A_2(x_0,\cdot)(\ux) \ = \ b_{-1}(\cdot) \ast A_2(x_0,\cdot)(\ux) \ = \ B_2(x_0,\ux) \\[2mm]
\mathcal{H} \left [ B_2(x_0,\ux) \right ] & = & \mathcal{H}^2 \left [ A_2(x_0,\ux) \right ] \ = \ b_{-1}(\cdot) \ast B_2(x_0,\cdot)(\ux) \ = \ A_2(x_0,\ux) 
\end{array} \right .
$$
while, trivially,
$$
\left \{ \begin{array}{rcl}
a_{-1}(\cdot) \ast A_2(x_0,\cdot)(\ux) &  = & A_2(x_0,\ux) \\[2mm]
a_{-1}(\cdot) \ast B_2(x_0,\cdot)(\ux) &  = & B_2(x_0,\ux)
\end{array} \right .
$$
Their distributional limits for $x_0 \rightarrow 0+$ are given by
$$
\left \{ \begin{array}{rcl}
a_2(\ux) & = & \lim_{x_0 \rightarrow 0+} A_2(x_0,\ux) \ = \ a_0 \ast a_1(\ux) \ = \ b_0 \ast b_1(\ux) \\[2mm]
b_2(\ux) & = & \lim_{x_0 \rightarrow 0+} B_2(x_0,\ux) \ = \ a_0 \ast b_1(\ux) \ = \ b_0 \ast a_1(\ux) 
\end{array} \right .
$$
which may be calculated explicitly to be
$$
a_2(\ux) = \left ( - \frac{2}{m-1} \frac{1}{\sigma_{m+1}} T^\ast_{-m+1} \right ) \ast \left ( \frac{1}{\pi} \frac{1}{\sigma_m} \frac{1}{m-2} T^\ast_{-m+2} \right ) = - \frac{2}{(m-1)(m-3)} \frac{1}{\sigma_{m+1}} \frac{1}{|\ux|^{m-3}}
$$
or
$$
a_2(\ux) = \left ( \frac{1}{\pi} \frac{1}{\sigma_{m}} U^\ast_{-m+1} \right ) \ast \left ( - \frac{1}{\pi} \frac{1}{\sigma_{m+1}} \frac{1}{m-1} U^\ast_{-m+2} \right ) = - \frac{1}{\pi} \frac{1}{(m-1)(m-3)} \frac{1}{\sigma_{m+1}} \, T^\ast_{-m+3}
$$
and
$$
b_2(\ux) = \left ( - \frac{2}{m-1} \frac{1}{\sigma_{m+1}} T^\ast_{-m+1} \right ) \ast \left ( - \frac{1}{\pi} \frac{1}{\sigma_{m+1}} \frac{1}{m-1} U^\ast_{-m+2} \right ) =  \frac{1}{2} \frac{1}{\sigma_{m}} \frac{1}{m-2} \frac{\ux}{|\ux|^{m-2}}
$$
or
$$
b_2(\ux) = \left ( \frac{1}{\pi} \frac{1}{\sigma_{m}} U^\ast_{-m+1} \right ) \ast \left ( \frac{1}{\pi} \frac{1}{\sigma_{m}} \frac{1}{m-2} T^\ast_{-m+2} \right ) =  \frac{1}{2 \pi^2} \frac{1}{\sigma_{m}} \frac{1}{m-2} U^\ast_{-m+3}
$$
Putting $c_2(\ux) = \frac{1}{2} a_2(\ux) + \frac{1}{2} \overline{e_0} b_2(\ux)$ we can prove the following properties of those distributional boundary values.
\begin{lemma}
\label{lem56}
\rule{0mm}{0mm}
\begin{itemize}
\item[(i)] $-\pux a_2 = b_1$, $- \pux b_2 = a_1$, $- \overline{e_0} \pux c_2 = c_1$
\item[(ii)] $\mathcal{H} \left [ a_2 \right ] = b_2$, $\mathcal{H} \left [ b_2 \right ] = a_2$, $\overline{e_0} \mathcal{H} \left [ c_2 \right ] = c_2$
\item[(iii)] $c_{-1} \ast a_2 = c_{-1} \ast b_2 = c_{-1} \ast c_2 = c_2$
\item[(iv)] $c_0 \ast a_1 = a_1 \ast c_0 = c_2$, $c_0 \ast \overline{e_0} b_1 = \overline{e_0} b_1 \ast c_0 = c_2$ 
\end{itemize}
\end{lemma}

\pf
For (i) we obtain
\begin{eqnarray*}
- \pux a_2(\ux) & = & - \pux \left ( - \frac{1}{\pi} \frac{1}{\sigma_{m+1}} \frac{1}{m-1} \frac{1}{m-3} \, T^\ast_{-m+3} \right ) \ = \ \frac{1}{\pi} \frac{1}{\sigma_{m+1}} \frac{1}{m-1} \frac{1}{m-3} (-m+3) \, U^\ast_{-m+2} \\
& = & - \frac{1}{\pi} \frac{1}{m-1} \frac{1}{\sigma_{m+1}} \, U^\ast_{-m+2} \ = \ b_1(\ux) \\
- \pux b_2(\ux) & = & - \pux \left ( - \frac{1}{2\pi^2} \frac{1}{\sigma_{m}} \frac{1}{m-2} \, U^\ast_{-m+3} \right ) \ = \ - \frac{1}{2\pi^2} \frac{1}{\sigma_{m}} \frac{1}{m-2} (-2\pi) \, T^\ast_{-m+2} \\
& = & \frac{1}{\pi} \frac{1}{m-2} \frac{1}{\sigma_{m}} \, T^\ast_{-m+2} \ = \ a_1(\ux) \\
\left ( \overline{e_0} \pux \right ) c_2(\ux) & = & \left ( \overline{e_0} \pux \right ) \left ( \frac{1}{2} a_2 + \frac{1}{2} \overline{e_0} b_2 \right ) = \frac{1}{2} \overline{e_0} b_1 + \frac{1}{2} a_1 \ = \ c_1(\ux)
\end{eqnarray*}
while for (ii)
\begin{eqnarray*}
\mathcal{H} \left [ a_2 \right ] & = & \mathcal{H} \left [ a_0 \ast a_1 \right ] \ = \ \mathcal{H} \left [ a_0 \right ] \ast a_1 \ = \ b_0 \ast a_1 \ = \ b_2 \\[1mm]
\mathcal{H} \left [ b_2 \right ] & = & \mathcal{H}^2 \left [ a_2 \right ] \ = \ a_2 \\
\overline{e_0} \mathcal{H} \left [ c_2 \right ] & = & \overline{e_0} \mathcal{H} \left [ \frac{1}{2} a_2 + \frac{1}{2} \overline{e_0} b_2 \right ] \ = \ \frac{1}{2} \overline{e_0} b_2 + \frac{1}{2} a_2 \ = \ c_2
\end{eqnarray*}
Statements (iii) and (iv) follow by direct computation.
\qed

Making use of the distributional boundary values $a_1(\ux)$, $b_1(\ux)$, $a_2(\ux)$ and $b_2(\ux)$ we may now prove by direct computation the following equivalent expressions for the conjugate harmonic potentials $A_2(x_0,\ux)$ and $B_2(x_0,\ux)$.
\begin{lemma}
\label{lem57}
One has, convolutions being taken in the variable $\ux \in \mR^{m}$:
\begin{itemize}
\item[(i)] $A_2 = a_1 \ast A_0 = b_1 \ast B_0 = a_2 \ast A_{-1} = b_2 \ast B_{-1}$
\item[(ii)] $B_2 = a_1 \ast B_0 = b_1 \ast A_0 = a_2 \ast B_{-1} = b_2 \ast A_{-1}$
\end{itemize}
\end{lemma}
Also note the commutative schemes
$$
\begin{array}{ccc}
A_2(x_0,\ux) & \xrightarrow{\hspace*{2mm} - \pux \hspace*{2mm}} & B_{1}(x_0,\ux) \\[2mm]
\hspace*{-8mm} ^{\hspace*{0.1mm}_{x_0 \rightarrow 0+}} \downarrow & & \downarrow \\
a_{2}(\ux) &  \xrightarrow{\hspace*{1mm} - \pux \hspace*{1mm}} & b_{1}(\ux) 
\end{array}
\quad \mbox{and} \qquad 
\begin{array}{ccc}
B_{2}(x_0,\ux) &  \xrightarrow{\hspace*{1mm} - \pux \hspace*{1mm}} & A_{1}(x_0,\ux)  \\[2mm]
\hspace*{-7mm} ^{\hspace*{0.1mm}_{x_0 \rightarrow 0+}} \downarrow & & \downarrow \\
b_{2}(\ux) &  \xrightarrow{\hspace*{1mm} - \pux \hspace*{1mm}} & a_{1}(\ux) 
\end{array}
$$
We also have explicitly determined the conjugate harmonic potentials $A_2(x_0,\ux)$ (for $m>3$) and $B_2(x_0,\ux)$ (for $m>2)$:
\begin{eqnarray*}
A_2(x_0,\ux) & = & \frac{2}{m-1} \frac{1}{\sigma_{m+1}} \frac{x_0}{|\ux|^{m-2}} \, F_{m-2} \left ( \frac{|\ux|}{x_0} \right ) - \frac{2}{m-1} \frac{1}{m-3} \frac{1}{\sigma_{m+1}} \frac{1}{|x|^{m-3}} \\
B_2(x_0,\ux) & = & \frac{1}{\sigma_{m+1}} \frac{\ux |x|^2}{|\ux|^{m}} \, F_{m} \left ( \frac{|\ux|}{x_0} \right ) - \frac{m-3}{m-1} \frac{1}{\sigma_{m+1}} \frac{\ux}{|\ux|^{m-2}} \,   F_{m-2} \left ( \frac{|\ux|}{x_0} \right )
\end{eqnarray*}

\subsection{The general case}

Inspired by the properties of the harmonic potentials $A_1(x_0,\ux)$, $B_1(x_0,\ux)$, $A_2(x_0,\ux)$ and $B_2(x_0,\ux)$, we define recursively, for general $k=1,2,3,\ldots$ , the following functions in $\mR^{m+1}_+$, the convolutions being taken in the variable $\ux \in \mR^{m}$:
\begin{eqnarray*}
A_k(x_0,\ux) & = & a_0 \ast A_{k-1} \ = \ a_1 \ast A_{k-2} \ = \ \ldots \ = \ a_{k-1} \ast A_0 \\
& = & b_0 \ast B_{k-1} \ = \ b_1 \ast B_{k-2} \ = \ \ldots \ = \ b_{k-1} \ast B_0 \\
B_k(x_0,\ux) & = & a_0 \ast B_{k-1} \ = \ a_1 \ast B_{k-2} \ = \ \ldots \ = \ a_{k-1} \ast B_0 \\
& = & b_0 \ast A_{k-1} \ = \ b_1 \ast A_{k-2} \ = \ \ldots \ = \ b_{k-1} \ast A_0 \\
\end{eqnarray*}
and
$$
C_k(x_0,\ux) = \frac{1}{2} A_k(x_0,\ux) + \frac{1}{2} \overline{e_0} B_k(x_0,\ux)
$$
Note that for $k=1,2$ we indeed recover the harmonic potentials studied in the previous subsections. In a similar way as above, it is now shown that $A_k(x_0,\ux)$, $B_k(x_0,\ux)$ and $C_k(x_0,\ux)$ satisfy the following equations:
\begin{itemize}
\item[(i)] $\p_{x_0} A_k = A_{k-1}$ \\[1mm]
$- \pux A_k = B_{k-1}$ \\[1mm]
$\overline{D} A_k = \frac{1}{2} \left ( \p_{x_0} - \overline{e_0} \pux \right ) A_k = \frac{1}{2} A_{k-1} + \frac{1}{2} \overline{e_0} B_{k-1} = C_{k-1}$
\item[(ii)] $\p_{x_0} B_k = B_{k-1}$ \\[1mm]
$- \pux B_k = A_{k-1}$ \\[1mm]
$\overline{D} \left ( \overline{e_0} B_k \right ) = \frac{1}{2} \left ( \p_{x_0} - \overline{e_0} \pux \right ) \overline{e_0} B_k = \frac{1}{2} \overline{e_0} B_{k-1} + \frac{1}{2} A_{k-1} = C_{k-1}$
\item[(iii)] $D C_k = \frac{1}{2} \left ( \p_{x_0} + \overline{e_0} \pux \right ) \left ( \frac{1}{2} A_k + \frac{1}{2} \overline{e_0} B_k \right ) = 0$
\item[(iv)] $\overline{D} C_k = \overline{D} \left ( \frac{1}{2} A_k + \frac{1}{2} \overline{e_0} B_k \right ) = C_{k-1}$
\end{itemize}
which clearly show that $A_k(x_0,\ux)$ and $B_k(x_0,\ux)$ are conjugate harmonic potentials of $C_{k-1}(x_0,\ux)$ in $\mR^{m+1}_+$, while $C_k(x_0,\ux)$ is a monogenic potential of the same $C_{k-1}(x_0,\ux)$ in $\mR^{m+1}_+$. Their distributional boundary values for $x_0 \rightarrow 0+$ are given by the recurrence relations
\begin{eqnarray*}
a_k(\ux) & = & a_0 \ast a_{k-1}\ = \ a_1 \ast a_{k-2} \ = \ \ldots \ = \  a_{k-1} \ast a_0 \\
& = & b_0 \ast b_{k-1}\ = \ b_1 \ast b_{k-2} \ = \ \ldots \ = \  b_{k-1} \ast b_0 \\
b_k(\ux) & = & a_0 \ast b_{k-1}\ = \ a_1 \ast b_{k-2} \ = \ \ldots \ = \  a_{k-1} \ast b_0 \\
& = & b_0 \ast a_{k-1}\ = \ b_1 \ast a_{k-2} \ = \ \ldots \ = \  b_{k-1} \ast a_0 
\end{eqnarray*}
for which the following explicit formulae may be deduced:
$$
\left \{ \begin{array}{rcl}
a_{2j} & = & -\displaystyle\frac{1}{2^{j-1}} \displaystyle\frac{1}{\pi^j} \displaystyle\frac{1}{(m-1)(m-3) \ldots (m-2j-1)} \displaystyle\frac{1}{\sigma_{m+1}} \, T^\ast_{-m+2j+1} \\[5mm]
a_{2j-1} & = & \phantom{-} \displaystyle\frac{1}{2^{j-1}} \displaystyle\frac{1}{\pi^j} \displaystyle\frac{1}{(m-2)(m-4) \ldots (m-2j)} \displaystyle\frac{1}{\sigma_{m}} \, T^\ast_{-m+2j} 
\end{array} \right .
$$
$$
\left \{ \begin{array}{rcl}
b_{2j} & = & \phantom{-} \displaystyle\frac{1}{2^{j}} \displaystyle\frac{1}{\pi^{j+1}} \displaystyle\frac{1}{(m-2)(m-4) \ldots (m-2j)} \displaystyle\frac{1}{\sigma_{m}} \, U^\ast_{-m+2j+1} \\[5mm]
b_{2j-1} & = & - \displaystyle\frac{1}{2^{j-1}} \displaystyle\frac{1}{\pi^j} \displaystyle\frac{1}{(m-1)(m-3) \ldots (m-2j+1)} \displaystyle\frac{1}{\sigma_{m+1}} \, U^\ast_{-m+2j} 
\end{array} \right .
$$
These distributional limits show the following, by now traditional, properties.
\begin{lemma}
\label{lem58}
One has for $k=1,2,\ldots$:
\begin{itemize}
\item[(i)] $- \pux a_k = b_{k-1}$
\item[(ii)] $- \pux b_k = a_{k-1}$
\item[(iii)] $\mathcal{H} \left [ a_k \right ] = b_{-1} \ast a_k = b_k$
\item[(iv)] $\mathcal{H} \left [ b_k \right ] = b_{-1} \ast b_k = a_k$
\end{itemize}
\end{lemma}

\pf
Follows by direct computation using the derivation and convolution formulae for the $T^\ast$-- an $U^\ast$--distributions.
\qed


\section{Conclusion}


While constructing a higher dimensional analogue in upper half--space $\mR^{m+1}_+$ of the function $\ln{z}$ in the upper half of the complex plane, preserving its fundamental property of being a holomorphic potential of the Cauchy kernel $\frac{1}{z}$, it became clear that this monogenic logarithmic function  is but one of a double sequence of such kind of potentials, just as $\ln{z}$ is the central element in the double sequence of holomorphic primitives:
$$
\frac{1}{k!} z^k \left[ \ln z - ( 1 + \frac{1}{2} + \ldots + \frac{1}{k}) \right] \rightarrow \ldots \rightarrow z ( \ln z - 1) \rightarrow \ln z 
\stackrel{\frac{d}{dz}}{\longrightarrow} \frac{1}{z} \rightarrow - \frac{1}{z^2} \rightarrow \ldots \rightarrow (-1)^{k-1} \frac{(k-1)!}{z^k}
$$
The sequence of monogenic potentials corresponding to the negative integer powers of  $z$, which we called {\em downstream potentials}, were rather easily constructed via differentiation with the conjugate generalized Cauchy--Riemann operator $\overline{D}$. The explicit construction of the monogenic potentials corresponding to the logarithmic functions in $\mC_+$, which we termed {\em upstream potentials}, requires tedious calculations involving primitivation with respect to $\overline{D}$, and up to now we have executed three inductive steps. A general expression for these upstream potentials is lacking, but their properties are known since they arise as convolutions of adjacent potentials with their distributional boundary values in $\mR^m$. Also with an eye on possible applications, the upstream potentials will be further calculated in the lower dimensional cases where $m=2,3$, and it is hoped for that a general formula, mimicking the one in the complex plane, will appear.\\[-2mm]

The above mentioned distributional boundary values are really fundamental, since not only they are used in the definition of the potentials, but also uniquely determine the conjugate harmonic potentials obtained by primitivation, thanks to the simple, but crucial, fact that a monogenic function in $\mR^{m+1}_+$ vanishing at the boundary $\mR^{m}$ indeed is zero. For those distributional boundary values we have established a general formula, showing that they all fit into two families of distributions in $\mR^{m}$, one scalar--valued, the second one Clifford vector--valued. In some particular cases they have been identified as fundamental solutions of the Dirac operator or the Laplace operator, or as convolution kernels for some pseudodifferential operators related to both these operators. The forthcoming paper \cite{bdbds} will treat this remarkable relationship between the distributional boundary values of the harmonic potentials and specific integer and half--integer powers of the Dirac and Laplace operators.


\end{document}